\documentclass[10pt,a4paper]{article}
\usepackage{amsfonts,amsmath,mathrsfs,amssymb}
\usepackage{tikz}
\usepackage{verbatim}
\usepackage{array}
\newtheorem{theorem}{Theorem}[section]
\newtheorem{lemma}[theorem]{Lemma}

\newtheorem{defin}{Definition}[section]
\newtheorem{examp}{Example}[section]
\newtheorem{remark}{Remark}[section]
\newtheorem{cor}[theorem]{Corollary}

\begin{document}

\title {Symmetry as a sufficient condition for a finite flex}

\author{Bernd Schulze\footnote{Supported by the DFG Research Unit 565 `Polyhedral Surfaces'.}\\ Institute of Mathematics, MA 6-2\\ TU Berlin\\ D-10623 Berlin, Germany
}

\maketitle

\begin{abstract}

We show that if the joints of a bar and joint framework $(G,p)$ are positioned as `generically' as possible subject to given symmetry constraints and $(G,p)$ possesses a `fully-symmetric' infinitesimal flex (i.e., the velocity vectors of the infinitesimal flex remain unaltered under all symmetry operations of $(G,p)$), then $(G,p)$ also possesses a finite flex which preserves the symmetry of $(G,p)$ throughout the path. This and other related results are obtained by symmetrizing techniques described by L. Asimov and B. Roth in their paper `The Rigidity Of Graphs' from 1978 and by using the fact that the rigidity matrix of a symmetric framework can be transformed into a block-diagonalized form by means of group representation theory. The finite flexes that can be detected with these symmetry-based methods can in general not be found with the analogous non-symmetric methods.

\end{abstract}

\section{Introduction}

A bar and joint framework 
is said to be `rigid' if, loosely speaking, it cannot be deformed continuously into another non-congruent framework while keeping the lengths of all bars fixed. Otherwise, the framework is said to be `flexible' \cite{asiroth, gluck, graver, W2}. The identification of  flexes (sometimes also called `finite flexes' or `mechanisms')  in frameworks is, in general, a very difficult problem. For frameworks whose joints lie in `generic' positions, however,  L. Asimov and B. Roth showed in 1978 that `rigidity' is equivalent to `infinitesimal rigidity' (see \cite{asiroth}). So, for `generic' frameworks, finite flexes can easily be detected, since the infinitesimal rigidity properties of a framework are completely described by the linear algebra of its rigidity matrix \cite{gss, W1, W2}.\\\indent
While the result of L. Asimov and B. Roth applies to `almost all' realizations of a given graph, it does, in general, not apply to frameworks that possess non-trivial symmetries, since the joints of a symmetric framework are typically forced to lie in special `non-generic' positions \cite{BS4, BS1}.
\\\indent In  this paper, we establish symmetric versions of the theorem of Asimov and Roth which allow us to detect finite flexes in symmetric frameworks that are as generic as possible subject to the given symmetry constraints. These results are based on the fact that the rigidity matrix of a symmetric framework $(G,p)$ can be transformed into a block-diagonalized form using techniques from group representation theory  \cite{KG2, BS2}. In this block-diagonalization of the rigidity matrix of $(G,p)$, each submatrix block corresponds to an irreducible representation of the point group $S$ of $(G,p)$, so that the (infinitesimal) rigidity analysis of $(G,p)$ can be broken up into independent subproblems, where each subproblem considers the relationship between external forces on the joints and resulting internal distortions in the bars of $(G,p)$ that share certain symmetry properties. In particular, the submatrix block that corresponds to the trivial irreducible representation of $S$ describes the relationship between external displacement vectors on the joints and internal distortion vectors in the bars of $(G,p)$ that are `fully-symmetric', i.e., unchanged under all symmetry operations in $S$.
By adapting the techniques described by  L. Asimov and B. Roth in \cite{asiroth} and applying them to the submatrix block corresponding to the trivial irreducible representation of $S$, rather than the entire rigidity matrix, we show that the existence of a `fully-symmetric' infinitesimal flex of $(G,p)$ also implies the existence of a finite flex of $(G,p)$ which preserves the symmetry of $(G,p)$ throughout the path, provided that $(G,p)$ is as generic as possible subject to the given symmetry constraints.
As a corollary of this result, one obtains the Proposition 1 stated (but not proven) in \cite{FG4} (see also \cite{KG1}).
By considering submatrix blocks that correspond to other irreducible representations of $S$, finite flexes of $(G,p)$ that preserve some, but not all of the symmetries of $(G,p)$ can also be detected.\\\indent
In order to apply our symmetry-adapted versions of the theorem of Asimov and Roth to a given framework $(G,p)$, we need to  detect infinitesimal flexes that possess certain symmetry properties (in particular, `fully-symmetric' infinitesimal flexes) in $(G,p)$. In many cases, this can be done with very little computational effort by means of the Fowler-Guest symmetry-extended version of Maxwell's rule described in \cite{FGsymmax} (see also \cite{cfgsw, BS2, BS4}). An alternate way of finding `fully-symmetric' infinitesimal flexes, as well as `fully-symmetric' self-stresses, in symmetric frameworks will be presented in \cite{SchW}.\\\indent
We note that there exist a number of famous and interesting examples of symmetric frameworks that can be shown to be flexible with our symmetry-based methods. These include the frameworks examined in \cite{bottema, bricard, consusp, FG4, BS4, tarnaiely}, for example. For each of these classes of frameworks, our new approach provides a much simpler proof for the existence of a finite symmetry-preserving flex than previous methods. In the final section of this paper, we demonstrate the efficiency of our results by showing the flexibility of two of the three types of `Bricard octhedra' \cite{bricard, Stachel}. It is shown in \cite{SchW} that new finite flexes can also be detected with our methods.

\section{Rigidity theoretic definitions and preliminaries}
\label{sec:rig}

\subsection{Rigidity}
\label{subsec:rig}

All graphs considered in this paper are finite graphs without loops or multiple edges. The \emph{vertex set} of a graph $G$ is denoted by $V(G)$ and the \emph{edge set} of $G$ is denoted by $E(G)$.\\\indent
A \emph{framework} in $\mathbb{R}^{d}$ is a pair $(G,p)$, where $G$ is a graph and $p: V(G)\to \mathbb{R}^{d}$ is a map with the property that $p(u) \neq p(v)$ for all $\{u,v\} \in E(G)$. We also say that $(G,p)$ is a $d$-dimensional \emph{realization} of the \emph{underlying graph} $G$ \cite{graver, gss, W1, W2}.
An ordered pair $\big(v,p(v)\big)$, where $v \in V(G)$, is a \emph{joint} of $(G,p)$, and an unordered pair $\big\{\big(u,p(u)\big),\big(v,p(v)\big)\big\}$ of joints, where $\{u,v\} \in E(G)$, is a \emph{bar} of $(G,p)$.\\\indent
Given the vertex set $V(G)=\{v_{1},\ldots,v_{n}\}$ of a graph $G$ and a map $p:V(G)\to \mathbb{R}^{d}$, it is often useful to identify $p$ with a vector in $\mathbb{R}^{dn}$ by using the order on $V(G)$. In this case we also refer to $p$ as a \emph{configuration} of $n$ points in $\mathbb{R}^{d}$. For $p(v_{i})$ we will frequently write $p_{i}$.
\\\indent For a fixed ordering of the edges of a graph $G$ with vertex set $V(G)=\{v_{1},v_{2},\ldots, v_{n}\}$, we define the \emph{edge function} $f_{G}:\mathbb{R}^{dn}\to \mathbb{R}^{|E(G)|}$  by
\begin{displaymath}
f_{G}\big(p_{1},\ldots,p_{n}\big)=\big(\ldots, \|p_{i}-p_{j}\|^2,\ldots\big)\textrm{, }
\end{displaymath}
where $\{v_{i},v_{j}\}\in E(G)$, $p_{i}\in \mathbb{R}^{d}$ for all $i=1,\ldots,n$, and $\| \cdot \|$ denotes the Euclidean norm in $\mathbb{R}^{d}$  \cite{asiroth, BS4}.\\\indent
If $(G,p)$ is a $d$-dimensional framework with $n$ vertices, then $f_{G}^{-1}\big(f_{G}(p)\big)$ is the set of all configurations $q$ of $n$ points in $\mathbb{R}^{d}$ with the property that corresponding bars of the frameworks $(G,p)$ and $(G,q)$ have the same length. In particular, we clearly have $f^{-1}_{K_{n}}\big(f_{K_{n}}(p)\big)\subseteq f_{G}^{-1}\big(f_{G}(p)\big)$, where $K_{n}$ is the complete graph on $V(G)$.

\begin{defin}\label{altdefmotions}
\emph{\cite{asiroth, RW1} Let $G$ be a graph with $n$ vertices and let $(G,p)$ be a framework in $\mathbb{R}^{d}$. A \emph{motion} of $(G,p)$ is a differentiable path $x:[0,1]\to \mathbb{R}^{dn}$ such that $x(0)=p$ and $x(t)\in f_{G}^{-1}\big(f_{G}(p)\big)$ for all $t\in [0,1]$.\\\indent A motion $x$ of $(G,p)$ is a \emph{rigid motion} if $x(t)\in f_{K_{n}}^{-1}\big(f_{K_{n}}(p)\big)$ for all $t\in [0,1]$ and a \emph{flex} of $(G,p)$ if $x(t)\notin f_{K_{n}}^{-1}\big(f_{K_{n}}(p)\big)$ for all $t\in (0,1]$.\\\indent  $(G,p)$ is \emph{rigid} if every motion of $(G,p)$ is a rigid motion, and \emph{flexible} otherwise.}
\end{defin}

\begin{figure}[htp]
\begin{center}
\begin{tikzpicture}[very thick,scale=1]
\tikzstyle{every node}=[circle, draw=black, fill=white, inner sep=0pt, minimum width=5pt];
\node (p1) at (0,0) {};
\node (p2) at (1.7,0) {};
\node (p3) at (0.5,1) {};
\draw(p1)--(p2);
\draw(p3)--(p2);
\draw(p1)--(p3);
\node [rectangle,draw=white, fill=white] (a) at (0.85,-0.6) {(a)};
\end{tikzpicture}
\hspace{0.5cm}
\begin{tikzpicture}[very thick,scale=1]
\tikzstyle{every node}=[circle, draw=black, fill=white, inner sep=0pt, minimum width=5pt];
\node (p1) at (0,0) {};
\node (p2) at (1.7,0) {};
\node (p3) at (1.7,1) {};
\node (p4) at (0,1) {};
\draw(p1)--(p2);
\draw(p3)--(p2);
\draw(p4)--(p3);
\draw(p4)--(p1);
\node [rectangle,draw=white, fill=white] (b) at (0.85,-0.6) {(b)};
\end{tikzpicture}
\hspace{0.5cm}
\begin{tikzpicture}[very thick,scale=1]
\tikzstyle{every node}=[circle, draw=black, fill=white, inner sep=0pt, minimum width=5pt];
\node (p1) at (0,0) {};
\node (p2) at (1.7,0) {};
\node (p3) at (1.7,1) {};
\node (p4) at (0,1) {};
\node [draw=black!40!white](p5) at (0.6,0.8) {};
\node [draw=black!40!white](p6) at (2.3,0.8) {};
\draw(p1)--(p2);
\draw(p3)--(p2);
\draw(p4)--(p3);
\draw(p4)--(p1);
\draw[black!40!white](p5)--(p6);
\draw[black!30!white](p1)--(p5);
\draw[black!30!white](p2)--(p6);
\draw[thick,->,black!80!white](0.1,1) arc (90:48:19pt);
\draw[thick,->,black!80!white](1.8,1) arc (90:48:19pt);
\node [rectangle,draw=white, fill=white] (c) at (0.85,-0.6) {(c)};
\end{tikzpicture}
\hspace{0.5cm}
\begin{tikzpicture}[very thick,scale=1]
\tikzstyle{every node}=[circle, draw=black, fill=white, inner sep=0pt, minimum width=5pt];
\node (p1) at (0,0) {};
\node (p2) at (1.7,0) {};
\node (p5) at (0.6,0.8) {};
\node (p6) at (2.3,0.8) {};
\draw(p5)--(p6);
\draw(p1)--(p2);
\draw(p1)--(p5);
\draw(p2)--(p6);
\node [rectangle,draw=white, fill=white] (d) at (0.85,-0.6) {(d)};
\end{tikzpicture}
\end{center}
\vspace{-0.3cm}
\caption{\emph{A rigid \emph{(b)} and a flexible \emph{(b)} framework in the plane. The flex shown in \emph{(c)} takes the framework in \emph{(b)} to the framework in \emph{(d)}. }}
\label{rigidexamples}
\end{figure}
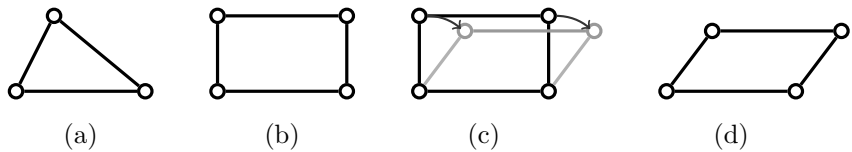

There also exist some alternate definitions of a flexible framework all of which are equivalent to Definition \ref{altdefmotions}, as shown in \cite{asiroth}.

\begin{theorem}
\label{alternaterigid}\cite{asiroth}
Let $(G,p)$ be a framework in $\mathbb{R}^{d}$ with $n$ vertices. The following are equivalent:
\begin{itemize}
\item[(i)] $(G,p)$ is flexible;
\item[(ii)] there exists a motion $x:[0,1]\to \mathbb{R}^{dn}$ of $(G,p)$ such that $x(t)\notin f_{K_{n}}^{-1}\big(f_{K_{n}}(p)\big)$ for some $t\in (0,1]$;
\item[(iii)] for every neighborhood $N_{p}$ of $p\in \mathbb{R}^{dn}$, we have $f_{K_{n}}^{-1}\big(f_{K_{n}}(p)\big)\cap N_{p}\varsubsetneqq f_{G}^{-1}\big(f_{G}(p)\big)\cap N_{p}$.
\end{itemize}
\end{theorem}

\begin{remark}
\emph{In Definition \ref{altdefmotions}, we may replace the term `differentiable path' by the terms `continuous path' or `analytic path'. The fact that all of these definitions are equivalent is a consequence of some basic results from algebraic geometry \cite{asiroth, RW1, W2}.}
\end{remark}

\subsection{Infinitesimal rigidity}
\label{subsec:infrig}

It is in general very difficult to determine whether a given framework $(G,p)$ is rigid or not since it requires solving a system of quadratic equations. It is therefore common to linearize this problem by considering the derivative of the edge function $f_G$ of $G$.

\begin{defin}\emph{\cite{graver, gss, W1, W2} Let $G$ be a graph with $V(G)=\{v_{1},v_{2},\ldots,v_{n}\}$ and let $p:V(G)\to \mathbb{R}^{d}$ be a map. The
\emph{rigidity matrix} of $(G,p)$ is the $|E(G)| \times dn$ matrix  \begin{displaymath} \mathbf{R}(G,p)=\left(
\begin{array} {ccccccccccc }
& & & & & \vdots & & & & & \\
0 & \ldots & 0 & p_{i}-p_{j}&0 &\ldots &0 & p_{j}-p_{i} &0 &\ldots &
0\\ & & & & & \vdots & & & & &\end{array}
\right)\textrm{,}\end{displaymath} that is, for each edge $\{v_{i},v_{j}\}\in E(G)$, $\mathbf{R}(G,p)$ has the row with
$(p_{i}-p_{j})_{1},\ldots,(p_{i}-p_{j})_{d}$ in the columns $d(i-1)+1,\ldots,di$, $(p_{j}-p_{i})_{1},\ldots,(p_{j}-p_{i})_{d}$ in
the columns $d(j-1)+1,\ldots,dj$, and $0$ elsewhere.\\\indent
Equivalently, $\mathbf{R}(G,p)=\frac{1}{2}df_{G}(p)$, where $df_{G}(p)$ denotes the Jacobian matrix of the edge function $f_{G}$ of $G$, evaluated at the point $p\in \mathbb{R}^{dn}$.}
\end{defin}

\begin{remark}
\emph{The rigidity matrix is defined for arbitrary pairs $(G,p)$, where $G$ is a graph and $p:V(G)\to \mathbb{R}^{d}$ is a map. If $(G,p)$ is not a framework, then there exists a pair of adjacent vertices of $G$ that are mapped to the same point in $\mathbb{R}^{d}$ under $p$ and every such edge of $G$ gives rise to a zero-row in $\mathbf{R}(G,p)$.}
\end{remark}

\begin{defin}\label{infdefmotions}
\emph{\cite{graver, gss, W1, W2} An \emph{infinitesimal motion} of a $d$-dimensional framework $(G,p)$ with $V(G)=\{v_{1},v_{2},\ldots,v_{n}\}$ is a function $u: V(G)\to \mathbb{R}^{d}$ such that
\begin{equation}
\label{infinmotioneq}
\big(p_i-p_j\big)\cdot \big(u_i-u_j\big)=0 \quad\textrm{ for all } \{v_i,v_j\} \in E(G)\textrm{,}\nonumber\end{equation} where $u_i=u(v_i)$.
Equivalently, $u$ is a vector in $\mathbb{R}^{dn}$ that lies in the kernel of $\mathbf{R}(G,p)$.
\\\indent An infinitesimal motion $u$ of $(G,p)$ is an \emph{infinitesimal rigid motion} if there exists a skew-symmetric matrix $S$ (a rotation) and a vector $t$ (a translation) such that $u_i=Sp_i+t$ for all $i=1,\ldots, n$; otherwise $u$ is an \emph{infinitesimal flex} of $(G,p)$.\\\indent $(G,p)$ is \emph{infinitesimally rigid} if every motion of $(G,p)$ is an infinitesimal rigid motion, and \emph{infinitesimally flexible} otherwise.}
\end{defin}

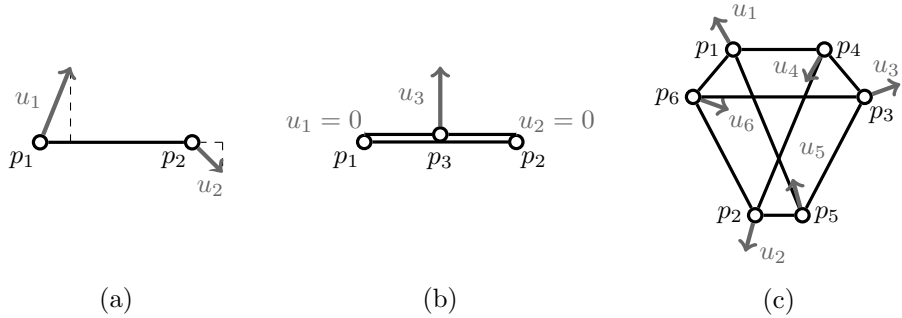
\begin{figure}[htp]
\begin{center}
\begin{tikzpicture}[very thick,scale=1]
\tikzstyle{every node}=[circle, draw=black, fill=white, inner sep=0pt, minimum width=5pt];
        \path (0,-0.2) node (p1) [label = below left: $p_{1}$] {} ;
        \path (2,-0.2) node (p2) [label = below left: $p_2$] {} ;
        \draw (p1)  --  (p2);
        \draw [dashed, thin] (0.4,-0.2) -- (0.4,0.8);
        \draw [dashed, thin] (2.4,-0.2) -- (2.4,-0.6);
        \draw [dashed, thin] (2.4,-0.2) -- (p2);
        \draw [->, black!60!white, ultra thick] (p1) -- node [draw=white, left=4pt] {$u_{1}$} (0.4,0.8);
        \draw [->, black!60!white, ultra thick] (p2) -- node [draw=white, below =4pt] {$u_2$} (2.4,-0.6);
        \node [draw=white, fill=white] (a) at (1,-2.3) {(a)};
        \end{tikzpicture}
        \hspace{0.5cm}
        \begin{tikzpicture}[very thick,scale=1]
\tikzstyle{every node}=[circle, draw=black, fill=white, inner sep=0pt, minimum width=5pt];
    \path (0,-0.2) node (p1) [label = below left: $p_1$] {} ;
    \path (2,-0.2) node (p2) [label = below right: $p_2$] {} ;
    \node (p3) at (1,-0.1) {};
    \node [draw=white, fill=white] (labelp3) at (1,-0.45) {$p_{3}$};
    \draw (p1) -- (p2);
    \draw (0,-0.1) -- (p3);
    \draw (p3) -- (2,-0.1);
    \draw [->, black!60!white, ultra thick] (p3) -- node [draw=white, left=4pt] {$u_{3}$} (1,0.8);
    \draw [ black!60!white, ultra thick] (p1) -- node [rectangle, draw=white, above left] {$u_{1}=0$} (p1);
    \draw [ black!60!white, ultra thick] (p2) -- node [rectangle, draw=white, above right] {$u_{2}=0$} (p2);
    \node [draw=white, fill=white] (b) at (1,-2.3) {(b)};
    \end{tikzpicture}
    \hspace{0.5cm}
        \begin{tikzpicture}[very thick,scale=1]
\tikzstyle{every node}=[circle, draw=black, fill=white, inner sep=0pt, minimum width=5pt];
    \path (160:1.2cm) node (p6) [label = left: $p_6$] {} ;
    \path (120:1.2cm) node (p1) [label = left: $p_1$] {} ;
    \path (255:1.2cm) node (p2) [label = left: $p_2$] {} ;
     \path (20:1.2cm) node (p3) [label = below right: $p_3$] {} ;
    \path (60:1.2cm) node (p4) [label = right: $p_4$] {} ;
    \path (285:1.2cm) node (p5) [label = right: $p_5$] {} ;
     \draw (p1) -- (p4);
     \draw [->, black!60!white, ultra thick] (p6) --  (160:0.7)node[draw=white, below right=0.3pt]{$u_{6}$};
     \draw (p1) -- (p5);
     \draw (p1) -- (p6);
     \draw (p2) -- (p4);
     \draw (p2) -- (p5);
     \draw (p2) -- (p6);
     \draw (p3) -- (p4);
          \draw (p3) -- (p6);
     \draw [->, black!60!white, ultra thick] (p1) -- node[rectangle, draw=white, above right=3pt] {$u_{1}$}(120:1.7);
     \draw [->, black!60!white, ultra thick] (p2) -- node[rectangle, draw=white, below right=3pt] {$u_{2}$}(255:1.7);
     \draw [->, black!60!white, ultra thick] (p3) -- node[rectangle, draw=white, above=4pt] {$u_{3}$}(20:1.7);
     \draw [->, black!60!white, ultra thick] (p4)  -- node[rectangle, draw=white, left=4pt] {$u_{4}$} node[rectangle, draw=white, below=27pt] {$u_{5}$}(60:0.7);
     \draw (p3) -- (p5);
     \draw [->, black!60!white, ultra thick] (p5) -- (285:0.7);
     \node [draw=white, fill=white] (c) at (0,-2.3) {(c)};
        \end{tikzpicture}
\end{center}
\vspace{-0.3cm}
\caption{\emph{The arrows indicate the non-zero displacement vectors of an infinitesimal rigid motion \emph{(a)} and infinitesimal flexes \emph{(b, c)} of frameworks in $\mathbb{R}^2$.}}
\label{inmo}
\end{figure}

Physically, an infinitesimal motion of $(G,p)$ is a set of initial velocity vectors, one at each joint, that preserve the lengths of the bars of $(G,p)$ at first order (see also Figure \ref{inmo}).\\\indent
It is well known that the infinitesimal rigid motions arising from $d$ translations and $\binom{d}{2}$ rotations of $\mathbb{R}^{d}$ form a basis for the space of infinitesimal rigid motions of $(G,p)$, provided that the points $p_{1},\ldots ,p_{n}$ span an affine subspace of $\mathbb{R}^{d}$ of dimension at least $d-1$ \cite{gss, W1}. Thus, for such a framework $(G,p)$, we have $\textrm{nullity }\big(\mathbf{R}(G,p)\big)\geq d+\binom{d}{2}=\binom{d+1}{2}$ and $(G,p)$ is infinitesimally rigid if and only if $\textrm{nullity } \big(\mathbf{R}(G,p)\big)=\binom{d+1}{2}$ or equivalently, $\textrm{rank }\big(\mathbf{R}(G,p)\big)=d |V(G)| - \binom{d+1}{2}$.

\begin{theorem}
\label{infinrigaff}\cite{asiroth, gluck}
A framework $(G,p)$ in $\mathbb{R}^d$ is infinitesimally rigid if and only if either $\textrm{rank }\big(\mathbf{R}(G,p)\big)=d |V(G)| - \binom{d+1}{2}$ or $G$ is a complete graph $K_n$ and the points $p(v)$, $v\in V(G)$, are affinely independent.
\end{theorem}

The following theorem gives the main connection between rigidity and infinitesimal rigidity. A proof of this result can be found in  \cite{asiroth}, \cite{gluck} or \cite {RW1}, for example.

\begin{theorem}
\label{infrigidtorigid}
If a framework $(G,p)$ is infinitesimally rigid, then $(G,p)$ is rigid.
\end{theorem}

\subsection{Basic `generic' results}
\label{subsec:basicgenres}

In 1978, L. Asimov and B. Roth showed that for `almost all' realizations of a given graph $G$, infinitesimal rigidity and rigidity are equivalent. We need the following definition.

\begin{defin}
\label{regular}
\emph{\cite{asiroth, BS4} Let $G$ be a graph with $n$ vertices and let $d\geq 1$ be an integer. A point $p\in \mathbb{R}^{dn}$ is said to be a \emph{regular point} of $G$ if there exists a neighborhood $N_{p}$ of $p$ in $\mathbb{R}^{dn}$ so that $\textrm{rank }\big(\mathbf{R}(G,p)\big)\geq\textrm{rank } \big(\mathbf{R}(G,q)\big)$ for all $q\in N_{p}$.\\\indent A framework $(G,p)$ is said to be \emph{regular} if $p$ is a regular point of $G$.}
\end{defin}

\begin{theorem}[Asimov, Roth, 1978]\label{asireg}\cite{asiroth}
Let $G$ be a graph with $n$ vertices and let $(G,p)$ be a $d$-dimensional framework. If $p\in \mathbb{R}^{dn}$ is a regular point of $G$, then $(G,p)$ is infinitesimally rigid if and only if $(G,p)$ is rigid.
\end{theorem}

Note that the set of all regular points of a graph $G$ in $\mathbb{R}^{dn}$ forms a dense open subset of $\mathbb{R}^{dn}$. Moreover, all regular realizations of $G$ share the same infinitesimal (and, by Theorem \ref{asireg}, also finite) rigidity properties. Regular frameworks are therefore sometimes also referred to as `generic' frameworks \cite{gss, RW1}. 
However, since in combinatorial (or generic) rigidity it is often useful to have a notion of `generic' that is invariant under addition or deletion of edges in $G$, a `generic' framework is frequently also defined as follows.


\begin{defin}
\label{generic}
\emph{\cite{graver, gss} Let $G$ be a graph with $n$ vertices, $d\geq 1$ be an integer, and $K_n$ be the complete graph on $V(G)$. Further, let $\mathbf{R}(n,d)$ be the matrix that is obtained from the rigidity matrix $\mathbf{R}(K_{n},p)$ of a $d$-dimensional realization $(K_{n},p)$ by replacing each $(p_{i})_{j}\in\mathbb{R}$ with the variable $(p'_{i})_{j}$. Then we say that $p\in \mathbb{R}^{dn}$ is \emph{generic} if the determinant of any submatrix of $\mathbf{R}(K_{n},p)$ is zero only if the determinant of the corresponding submatrix of $\mathbf{R}(n,d)$ is (identically) zero.\\\indent
The framework $(G,p)$ is said to be \emph{generic} if $p$ is generic.}
\end{defin}

Like the set of all regular points of $G$, the set of all generic points is also a dense open subset of $\mathbb{R}^{dn}$ \cite{gss}.
Moreover, since a generic framework is clearly also regular, we immediately obtain the following corollary of Theorem \ref{asireg}.

\begin{cor}\label{asigen}
If a framework $(G,p)$ is generic, then $(G,p)$ is infinitesimally rigid if and only if $(G,p)$ is rigid.
\end{cor}

An easy but often useful observation concerning generic frameworks is that if a framework $(G,p)$ in $\mathbb{R}^d$ is generic (in the sense of Definition \ref{generic}), then the joints of $(G,p)$ are in \emph{general position}, that is, for $1\leq m\leq d$, no $m+1$ joints of $(G,p)$ lie in an $(m-1)$-dimensional affine subspace of $\mathbb{R}^d$ \cite{gss}.\\\indent  For further information on generic frameworks, we refer the reader to \cite{graver, gss, Wgeneric}, for example.


\begin{defin}
\label{indep}
\emph{\cite{ W4, W1} A \emph{self-stress} of a framework $(G,p)$ is a vector $\omega\in \mathbb{R}^{|E(G)|}$ that satisfies $\mathbf{R}(G,p)^{T} \omega=0$. If $(G,p)$ has a non-zero self-stress, then $(G,p)$ is said to be \emph{dependent} (since in this case there exists a linear dependency among the row vectors of $\mathbf{R}(G,p)$). Otherwise, $(G,p)$ is said to be \emph{independent}.}
\end{defin}

An independent framework is clearly also regular. Therefore, the next result is also an immediate consequence of Theorem \ref{asireg}.

\begin{cor}\label{asiind}
If a framework $(G,p)$ is independent, then $(G,p)$ is infinitesimally rigid if and only if $(G,p)$ is rigid.
\end{cor}

We establish symmetric analogs to the Theorem of Asimov and Roth, as well as to Corollaries \ref{asigen} and \ref{asiind}, in Section \ref{sec:main}. These results will allow us to detect flexes in frameworks that are `generic' with respect to certain symmetry constraints.

\section{Symmetric frameworks}
\label{sec:sym}

\subsection{The set $\mathscr{R}_{(G,S,\Phi)}$}
\label{subsec:class}

\begin{figure}[htp]
\begin{center}
\begin{tikzpicture}[very thick,scale=1]
\tikzstyle{every node}=[circle, draw=black, fill=white, inner sep=0pt, minimum width=5pt];
    \path (0.3,-0.5) node (p5) [label = below left: $p_{5}$] {} ;
    \path (1.5,-1) node (p3) [label = below left: $p_{3}$] {} ;
    \path (2.7,-0.5) node (p6) [label = below right: $p_{6}$] {} ;
   \path (0.8,0.7) node (p1) [label = left: $p_{1}$] {} ;
   \path (2.2,0.7) node (p2) [label = right: $p_{2}$] {} ;
   \path (1.5,1.1) node (p4) [label = above left: $p_{4}$] {} ;
   \draw (p1) -- (p4);
     \draw (p1) -- (p5);
     \draw (p1) -- (p6);
     \draw (p2) -- (p4);
     \draw (p2) -- (p5);
     \draw (p2) -- (p6);
     \draw (p3) -- (p4);
     \draw (p3) -- (p5);
     \draw (p3) -- (p6);
   \draw [dashed, thin] (1.5,-1.7) -- (1.5,1.7);
      \node [draw=white, fill=white] (a) at (1.5,-2.2) {(a)};
    \end{tikzpicture}
    \hspace{2cm}
        \begin{tikzpicture}[very thick,scale=1]
\tikzstyle{every node}=[circle, draw=black, fill=white, inner sep=0pt, minimum width=5pt];
    \path (160:1.2cm) node (p6) [label = left: $p_6$] {} ;
    \path (120:1.2cm) node (p1) [label = above left: $p_1$] {} ;
    \path (255:1.2cm) node (p2) [label = below left: $p_2$] {} ;
     \path (20:1.2cm) node (p3) [label = right: $p_3$] {} ;
    \path (60:1.2cm) node (p4) [label = above right: $p_4$] {} ;
    \path (285:1.2cm) node (p5) [label = below right: $p_5$] {} ;
     \draw (p1) -- (p4);
     \draw (p1) -- (p5);
     \draw (p1) -- (p6);
     \draw (p2) -- (p4);
     \draw (p2) -- (p5);
     \draw (p2) -- (p6);
     \draw (p3) -- (p4);
     \draw (p3) -- (p5);
     \draw (p3) -- (p6);
     \draw [dashed, thin] (0,-1.7) -- (0,1.7);
      \node [draw=white, fill=white] (b) at (0,-2.2) {(b)};
        \end{tikzpicture}
\end{center}
\vspace{-0.3cm}
\caption{\emph{$2$-dimensional realizations of $K_{3,3}$ in $\mathscr{R}_{(K_{3,3},\mathcal{C}_s)}$ of different types: the framework in \emph{(a)} is of type
$\Phi_{a}$, where $\Phi_{a}: \mathcal{C}_{s} \to \textrm{Aut}(K_{3,3})$ is the homomorphism defined by $\Phi_{a}(s)=
(v_{1}\,v_{2})(v_{5}\,v_{6})(v_{3})(v_{4})$ and the framework in \emph{(b)} is of type
$\Phi_{b}$, where $\Phi_{b}: \mathcal{C}_{s} \to \textrm{Aut}(K_{3,3})$ is the homomorphism defined by $\Phi_{b}(s)=
(v_{1}\,v_{4})(v_{2}\,v_{5})(v_{3}\,v_{6})$.}}
\label{K33types}
\end{figure}
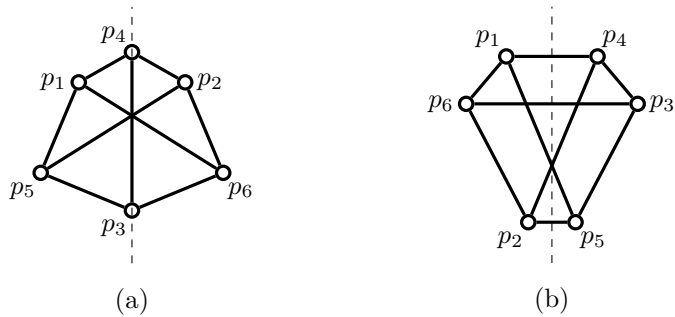

Recall that an \emph{automorphism} of a graph $G$ is a permutation $\alpha$ of $V(G)$ such that $\{u,v\}\in E(G)$ if and only if $\{\alpha(u),\alpha(v)\}\in E(G)$. The automorphisms of a graph $G$ form a group under composition which is denoted by $\textrm{Aut}(G)$.\\\indent
A \emph{symmetry operation} of a $d$-dimensional framework $(G,p)$ is an isometry $x$ of $\mathbb{R}^{d}$ such that for some $\alpha\in \textrm{Aut}(G)$, we have
$x\big(p(v)\big)=p\big(\alpha(v)\big)$ for all $v\in V(G)$ \cite{BS2, BS4, BS1}.\\\indent
The set of all symmetry operations of a framework $(G,p)$ forms a group under composition, called the \emph{point group} of $(G,p)$ \cite{bishop, Hall, BS4, BS1}. Since translating a framework does not change its rigidity properties, we may assume wlog that the point group of any framework in this paper is a \emph{symmetry group}, i.e., a subgroup of the orthogonal group $O(\mathbb{R}^{d})$.\\\indent
We use the Schoenflies notation for the symmetry operations and symmetry groups considered in this paper, as this is one of the standard notations in the literature about symmetric structures (see \cite{cfgsw, FGsymmax, FG4, KG1,  BS2, BS4, BS1, BS5}, for example). In this notation, the identity transformation is denoted by $Id$, a rotation  about a $(d-2)$-dimensional subspace of $\mathbb{R}^d$ by an angle of $\frac{2\pi}{m}$ is denoted by $C_m$, and a reflection in a $(d-1)$-dimensional subspace of $\mathbb{R}^d$ is denoted by $s$.\\\indent In this paper, we only consider three types of symmetry groups. In the Schoenflies notation, they are denoted by $\mathcal{C}_{s}$, $\mathcal{C}_{m}$, and $\mathcal{C}_{mv}$; $\mathcal{C}_{s}$ is a symmetry group consisting of the identity $Id$ and a single reflection $s$, $\mathcal{C}_{m}$ is a cyclic group generated by a rotation $C_m$, and $\mathcal{C}_{mv}$ is a dihedral group generated by a pair $\{C_m,s\}$. For further information about the Schoenflies notation we refer the reader to \cite{bishop, Hall, BS4}.
\\\indent
Given a symmetry group $S$ in dimension $d$ and a graph $G$, we let $\mathscr{R}_{(G,S)}$ denote the set of all $d$-dimensional realizations of $G$ whose point group is either equal to $S$ or contains $S$ as a subgroup \cite{BS2, BS4, BS1}. In other words, the set $\mathscr{R}_{(G,S)}$ consists of all realizations $(G,p)$ of $G$ for which there exists a map $\Phi:S\to \textrm{Aut}(G)$ so that
\begin{equation}\label{class} x\big(p(v)\big)=p\big(\Phi(x)(v)\big)\textrm{ for all } v\in V(G)\textrm{ and all } x\in S\textrm{.}\end{equation}
A framework $(G,p)\in \mathscr{R}_{(G,S)}$ satisfying the equations in (\ref{class}) for the map $\Phi:S\to \textrm{Aut}(G)$ is said to be \emph{of type $\Phi$}, and the set of all realizations in $\mathscr{R}_{(G,S)}$ which are of type $\Phi$ is denoted by $\mathscr{R}_{(G,S,\Phi)}$ (see again \cite{BS2, BS4, BS1} as well as Figure \ref{K33types}).

\begin{remark}
\emph{If the map $p$ of a framework $(G,p)\in \mathscr{R}_{(G,S)}$ is non-injective, then $(G,p)$ can possibly be of more than just one type and a given type may not be a homomorphism. However, if $p$ is injective, then $(G,p)\in \mathscr{R}_{(G,S)}$ is of a unique type $\Phi$ and $\Phi$ is necessarily also a homomorphism. See \cite{BS1} for further details.}
\end{remark}

\subsection{The notion of $(S,\Phi)$-generic}
\label{subsec:symgeneric}

Recall from equation (\ref{class}) in the previous section that for every framework $(G,p)$ in the set $\mathscr{R}_{(G,S,\Phi)}$  with $V(G)=\{v_1,\ldots,v_n\}$,  we have $x\big(p(v_{i})\big)=p\big(\Phi(x)(v_{i})\big)$ for all $i=1,2,\ldots, n$ and all $x\in S$.
Since every element of $S$ is an orthogonal linear transformation, we may identify each $x\in S$ with its corresponding orthogonal matrix $M_{x}$ that represents $x$ with respect to the canonical basis of $\mathbb{R}^{d}$. Therefore, for each $x\in S$, the equations in (\ref{class}) can be written as
\begin{displaymath} \mathbf{M^{(x)}} \left( \begin{array} {c} p_{1} \\ p_{2}
\\ \vdots \\ p_{n} \end{array} \right) =  \mathbf{P_{\Phi(x)}} \left( \begin{array} {c} p_{1} \\ p_{2}
\\ \vdots \\ p_{n} \end{array} \right) \textrm{,}\end{displaymath}
where
\begin{displaymath} \mathbf{M^{(x)}}=\left( \begin{array} {cccc}
M_{x} & 0 & \ldots & 0 \\ 0 & M_{x} & \ddots & \vdots\\
\vdots & \ddots & \ddots & 0\\ 0 & \ldots &
0& M_{x} \end{array} \right)\textrm{,}\end{displaymath}
and $\mathbf{P_{\Phi(x)}}$ is the $dn \times dn$ matrix which is
obtained from the permutation matrix corresponding to $\Phi(x)$
by replacing each $1$ by a $d \times d$ identity matrix and each
$0$ by a $d \times d$ zero matrix.\\\indent
We denote $L_{x,\Phi}=\textrm{ker }\big(\mathbf{M^{(x)}}- \mathbf{P_{\Phi(x)}}\big)$ and $U= \bigcap_{x \in S} L_{x,\Phi}$. Then $U$ is a linear subspace of $\mathbb{R}^{dn}$ which may be interpreted as the space of all those (possibly non-injective) configurations of $n$ points in $\mathbb{R}^{d}$ that possess the symmetry imposed by $S$ and $\Phi$ \cite{BS4, BS1}. In particular, if $(G,p)$ is a framework in $\mathscr{R}_{(G,S,\Phi)}$, then the configuration $p$ is an element of $U$.
Therefore, if we fix a basis $\mathscr{B}_{U}=\{u_{1},u_{2},\ldots,u_{k}\}$ of $U$, then every framework $(G,p) \in
\mathscr{R}_{(G,S,\Phi)}$ can be represented uniquely by the $k\times 1$ coordinate vector of $p$ relative to $\mathscr{B}_{U}$.\\\indent
We let $\mathbf{R}_{\mathscr{B}_{U}}(n,d)$ denote the matrix that is obtained from the `indeterminate' rigidity matrix $\mathbf{R}(n,d)$ (see Definition \ref{generic}) by introducing a $k$-tuple $(t'_{1},t'_{2},\ldots ,t'_{k})$ of variables and replacing the $dn$ variables $(p'_{i})_{j}$ of $\mathbf{R}(n,d)$ as follows.\\ \indent
For each $i=1,2, \ldots ,n$ and each $j=1,\ldots,d$, we replace the variable $(p'_{i})_{j}$ in $\mathbf{R}(n,d)$ by the linear combination
$t'_{1}(u_{1})_{i_{j}}+t'_{2}(u_{2})_{i_{j}}+ \ldots + t'_{k}(u_{k})_{i_{j}}$.

\begin{remark}
\label{symindetrigmatrem}
\emph{\cite{BS4, BS1} Let $(G,p)\in \mathscr{R}_{(G,S,\Phi)}$ and $\mathscr{B}_{U}=\{u_{1},u_{2},\ldots,u_{k}\}$ be a basis of
$U=\bigcap_{x \in S} L_{x,\Phi}$. Then \begin{displaymath}\left( \begin{array} {c} p_{1} \\ p_{2} \\ \vdots \\
p_{n}\end{array} \right)=t_{1}u_{1}+\ldots+ t_{k}u_{k}\textrm{, for some }t_{1},\ldots,t_{k}\in \mathbb{R}\textrm{.}\end{displaymath}
So, if for $i=1,\ldots,k$, the variable $t'_{i}$ in $\mathbf{R}_{\mathscr{B}_{U}}(n,d)$ is replaced by
$t_{i}$ then we obtain the rigidity matrix $\mathbf{R}(K_{n},p)$ of the framework $(K_{n},p)$.}
\end{remark}


The following symmetry-adapted notion of `generic' for the set $\mathscr{R}_{(G,S,\Phi)}$  was introduced in \cite{BS1} (see also \cite{BS4, BS5}).

\begin{defin}
\label{symgeneric}
\emph{\cite{BS4, BS1, BS5} Let $G$ be a graph with $V(G)=\{v_{1},v_{2},\ldots,v_{n}\}$, $K_{n}$ be the complete graph with $V(K_{n})=V(G)$, $S$ be a symmetry group in dimension $d$, $\Phi$ be a map from $S$ to $\textrm{Aut}(G)$, and $\mathscr{B}_{U}$ be a basis of $U=\bigcap_{x \in S} L_{x,\Phi}$.\\\indent A point $p\in \mathbb{R}^{dn}$ is called \emph{$(S,\Phi,\mathscr{B}_{U})$-generic} if the determinant of any submatrix of $\mathbf{R}(K_{n},p)$ is equal to zero only if the determinant of the
corresponding submatrix of $\mathbf{R}_{\mathscr{B}_{U}}(n,d)$ is (identically) zero.\\\indent The point $p$ is said to be \emph{$(S,\Phi)$-generic} if $p$ is $(S,\Phi,\mathscr{B}_{U})$-generic for some basis $\mathscr{B}_{U}$ of $U$.\\\indent
\indent A framework $(G,p) \in \mathscr{R}_{(G,S,\Phi)}$ is \emph{$(S,\Phi,\mathscr{B}_{U})$-generic} if $p$ is $(S,\Phi,\mathscr{B}_{U})$-generic, and $(G,p)$ is \emph{$(S,\Phi)$-generic} if $(G,p)$ is
$(S,\Phi,\mathscr{B}_{U})$-generic for some basis $\mathscr{B}_{U}$ of $U$.}
\end{defin}

\begin{remark}
 \emph{It is shown in \cite{BS1} that the definition of $(S,\Phi)$-generic is independent of the choice of the basis of $U$.}
\end{remark}

Intuitively, an $(S,\Phi)$-generic realization of a graph $G$ is obtained by placing the vertices of a set of representatives for the symmetry orbits $S v=\{\Phi(x)(v)|\, x\in S\}$ into `generic' positions. The positions for the remaining vertices of $G$ are then uniquely determined by the symmetry constraints imposed by $S$ and $\Phi$ (see \cite{BS4, BS1}, for further details).\\\indent As shown in \cite{BS1}, the set of $(S,\Phi)$-generic realizations of a graph $G$ is a dense open subset of the set $\mathscr{R}_{(G,S,\Phi)}$ and the infinitesimal rigidity properties are the same for all $(S,\Phi)$-generic realizations of $G$.
\\\indent
As an example, consider the realizations of $K_{3,3}$ with mirror symmetry depicted in Figure \ref{K33types}. All $(\mathcal{C}_s,\Phi_a)$-generic realizations of $K_{3,3}$  are infinitesimally rigid, whereas all realizations in $\mathscr{R}_{(K_{3,3},\mathcal{C}_s,\Phi_b)}$ are infinitesimally flexible (though rigid), since the joints of any realization in $\mathscr{R}_{(K_{3,3},\mathcal{C}_s,\Phi_b)}$ are forced to lie on a conic section \cite{BS1, W3}.

\subsection{The external and internal representation}
\label{subsec: blockdia}

Basic to our investigation of the rigidity and flexibility of symmetric frameworks in $\mathscr{R}_{(G,S,\Phi)}$ are the `external' and `internal' representation of the group $S$ which were first introduced in \cite{FGsymmax} and \cite{KG2} by means of an example (see also \cite{BS2,BS4}). We need mathematically rigorous definitions of these representations.


\begin{defin}\emph{ Let $G$ be a graph with $V(G)=\{v_{1},v_{2},\ldots,v_{n}\}$ and $E(G)=\{e_{1},e_{2},\ldots,e_{m}\}$, $S$ be a symmetry group in dimension $d$, and $\Phi:S\to \textrm{Aut}(G)$ be a homomorphism.\\\indent
The \emph{external representation} of $S$ (with respect to $G$ and $\Phi$) is the matrix representation $H_{e}:S\to GL(dn,\mathbb{R})$ that sends $x\in S$ to the matrix $H_{e}(x)$ which is obtained from the transpose of the $n\times n$ permutation matrix corresponding to $\Phi(x)$ (with respect to the enumeration $V(G)=\{v_{1},v_{2},\ldots,v_{n}\}$) by replacing each 1 with the orthogonal $d\times d$ matrix $M_{x}$ which represents $x$ with respect to the canonical basis of $\mathbb{R}^{d}$ and each 0 with a $d\times d$ zero-matrix.\\\indent We let $H'_{e}:S\to GL(\mathbb{R}^{dn})$ be the corresponding linear representation of $S$ that sends $x\in S$ to the automorphism $H'_{e}(x)$ which is represented by the matrix $H_{e}(x)$ with respect to the canonical basis of the $\mathbb{R}$-vector space $\mathbb{R}^{dn}$.\\\indent The \emph{internal representation} of $S$ (with respect to $G$ and $\Phi$) is the matrix representation $H_{i}:S\to GL(m,\mathbb{R})$ that sends $x\in S$ to the transpose of the permutation matrix corresponding to the permutation of $E(G)$ (with respect to the enumeration $E(G)=\{e_{1},e_{2},\ldots,e_{m}\}$) which is induced by $\Phi(x)$.\\\indent We let $H'_{i}:S\to GL(\mathbb{R}^{m})$ be the corresponding linear representation of $S$ that sends $x\in S$ to the automorphism $H'_{i}(x)$ which is represented by the matrix $H_{i}(x)$ with respect to the canonical basis of the $\mathbb{R}$-vector space $\mathbb{R}^{m}$.}
\end{defin}

\begin{figure}[htp]
\begin{center}
\begin{tikzpicture}[very thick,scale=1]
\tikzstyle{every node}=[circle, draw=black, fill=white, inner sep=0pt, minimum width=5pt];
        \path (0,0) node (p1) [label = below left: $p_1$] {} ;
        \path (2,0) node (p2) [label = below right: $p_2$] {} ;
        \path (1,2) node (p3) [label = left: $p_3$] {} ;
        \draw (p1)  -- node [draw=white, below left=3pt] {$e_{1}$} (p2);
        \draw (p2) -- node [draw=white, right=2pt] {$e_{3}$} (p3);
        \draw (p3) -- node [draw=white, left=2pt] {$e_{2}$} (p1);
        \draw [dashed, thin] (1,-0.6) -- (1,2.6);
        \end{tikzpicture}
\end{center}
\vspace{-0.3cm}
\caption{\emph{A framework $(K_{3},p)\in \mathscr{R}_{(K_{3},\mathcal{C}_{s},\Phi)}$.}}
\label{triang}
\end{figure}
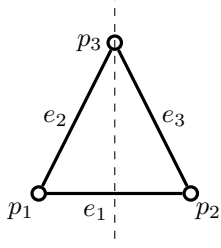

\begin{examp}
\label{triangexam}
\emph{To illustrate the previous definition, let $K_{3}$ be the complete graph with $V(K_{3})=\{v_{1},v_{2},v_{3}\}$ and $E(K_{3})=\{e_{1},e_{2},e_{3}\}$, where $e_{1}=\{v_{1},v_{2}\}$, $e_{2}=\{v_{1},v_{3}\}$ and $e_{3}=\{v_{2},v_{3}\}$. Further, let $\mathcal{C}_{s}=\{Id,s\}$ be the symmetry group in dimension 2 with \begin{displaymath} M_{Id}= \left( \begin{array} {rr} 1 & 0 \\ 0 & 1 \end{array} \right)\textrm{ and } M_{s}= \left( \begin{array} {rr} -1 & 0 \\ 0 & 1 \end{array} \right)\textrm{,}\end{displaymath} and let $\Phi:\mathcal{C}_{s}\to \textrm{Aut}(K_{3})$ be the homomorphism defined by $\Phi(s)=(v_{1}\,v_{2})(v_{3})$. (See also Figure \ref{triang}.) Then we have}
\begin{displaymath}
H_{e}(Id)= \left( \begin{array} {rr|rr|rr} 1 & 0 & 0 & 0 & 0 & 0\\ 0 & 1 & 0 & 0 & 0 & 0\\\hline 0 & 0 & 1 & 0 & 0 & 0\\0 & 0 & 0 & 1 & 0 & 0\\\hline 0 & 0 & 0 & 0 & 1 & 0\\0 & 0 & 0 & 0 & 0 & 1 \end{array} \right)\textrm{, } H_{e}(s)= \left( \begin{array} {rr|rr|rr} 0 & 0 & -1 & 0 & 0 & 0\\ 0 & 0 & 0 & 1 & 0 & 0\\\hline -1 & 0 & 0 & 0 & 0 & 0\\0 & 1 & 0 & 0 & 0 & 0\\\hline 0 & 0 & 0 & 0 & -1 & 0\\0 & 0 & 0 & 0 & 0 & 1\\\end{array} \right)\textrm{,}
\end{displaymath}
\begin{displaymath}
H_{i}(Id)= \left( \begin{array} {rrr} 1 & 0 & 0 \\ 0 & 1 & 0 \\ 0 & 0 & 1 \end{array} \right)\textrm{, } H_{i}(s)= \left( \begin{array} {rrr} 1 & 0 & 0\\ 0 & 0 & 1\\ 0 & 1 & 0\\\end{array} \right)\textrm{.}
\end{displaymath}
\end{examp}

For further examples, see \cite{KG2} or \cite{KG3}.\\\indent
From group representation theory we know that every finite group has, up to equivalency, only finitely many irreducible linear representations and that every linear representation of such a group can be written uniquely, up to equivalency of the direct summands, as a direct sum of the irreducible linear representations of this group \cite{liebeck, serre}. So, let $S$ have $r$ pairwise non-equivalent irreducible linear representations $I_{1},I_{2},\ldots, I_{r}$ and let
\begin{equation}
\label{irrrep}
H'_{e}= \lambda_{1}I_{1}\oplus\ldots\oplus \lambda_{r}I_{r} \textrm{, where } \lambda_{1},\ldots,\lambda_{r}\in \mathbb{N}\cup {\{0\}} \textrm{.}
\end{equation}
Then, for each $t=1,\ldots,r$, there exist $\lambda_{t}$ subspaces $\big(V_{e}^{(I_{t})}\big)_{1},\ldots, \big(V_{e}^{(I_{t})}\big)_{\lambda_{t}}$ of the $\mathbb{R}$-vector space $\mathbb{R}^{dn}$ which correspond to the $\lambda_{t}$ direct summands in (\ref{irrrep}), so that
\begin{equation}
\label{dirsumofvs}
\mathbb{R}^{dn}=V_{e}^{(I_{1})}\oplus \ldots \oplus V_{e}^{(I_{r})} \textrm{,}
\end{equation}
where
\begin{equation}
\label{dirsumofvs2}
V_{e}^{(I_{t})}= \big(V_{e}^{(I_{t})}\big)_{1}\oplus \ldots \oplus \big(V_{e}^{(I_{t})}\big)_{\lambda_{t}} \textrm{.}
\end{equation}
Let $\big(B_{e}^{(I_{t})}\big)_{1},\ldots, \big(B_{e}^{(I_{t})}\big)_{\lambda_{t}}$ be bases of the subspaces in (\ref{dirsumofvs2}). Then
\begin{equation}
\label{bases}
B_{e}^{(I_{t})}= \big(B_{e}^{(I_{t})}\big)_{1}\cup \ldots \cup \big(B_{e}^{(I_{t})}\big)_{\lambda_{t}}\nonumber
\end{equation}
is a basis of $V_{e}^{(I_{t})}$ and
\begin{equation}
\label{bases2}
B_{e}= B_{e}^{(I_{1})}\cup \ldots \cup B_{e}^{(I_{r})}
\end{equation}
is a basis of the $\mathbb{R}$-vector space $\mathbb{R}^{dn}$.
\\\indent
Similarly, let
\begin{equation}
\label{irrrepHi}
H'_{i}= \mu_{1}I_{1}\oplus\ldots\oplus \mu_{r}I_{r} \textrm{, where } \mu_{1},\ldots,\mu_{r}\in \mathbb{N}\cup {\{0\}} \textrm{.}
\end{equation}
For each $t=1,\ldots, r$, there exist $\mu_{t}$ subspaces $\big(V_{i}^{(I_{t})}\big)_{1},\ldots, \big(V_{i}^{(I_{t})}\big)_{\mu_{t}}$ of the $\mathbb{R}$-vector space $\mathbb{R}^{m}$ which correspond to the $\mu_{t}$ direct summands in (\ref{irrrepHi}), so that
\begin{equation}
\label{dirsumofvsHi}
\mathbb{R}^{m}=V_{i}^{(I_{1})}\oplus \ldots \oplus V_{i}^{(I_{r})} \textrm{,}
\end{equation}
where
\begin{equation}
\label{dirsumofvs2Hi}
V_{i}^{(I_{t})}= \big(V_{i}^{(I_{t})}\big)_{1}\oplus \ldots \oplus \big(V_{i}^{(I_{t})}\big)_{\mu_{t}} \textrm{.}
\end{equation}
Let $\big(B_{i}^{(I_{t})}\big)_{1},\ldots, \big(B_{i}^{(I_{t})}\big)_{\mu_{t}}$ be bases of the subspaces in (\ref{dirsumofvs2Hi}). Then
\begin{equation}
\label{basesHi}
B_{i}^{(I_{t})}= \big(B_{i}^{(I_{t})}\big)_{1}\cup \ldots \cup \big(B_{i}^{(I_{t})}\big)_{\mu_{t}}\nonumber
\end{equation}
is a basis of $V_{i}^{(I_{t})}$ and
\begin{equation}
\label{bases2Hi}
B_{i}= B_{i}^{(I_{1})}\cup \ldots \cup B_{i}^{(I_{r})}
\end{equation}
is a basis of the $\mathbb{R}$-vector space $\mathbb{R}^{m}$.

\begin{defin}\label{usym}
\emph{\cite{BS2} With the notation above, we say that a vector $v\in \mathbb{R}^{dn}$ is \emph{symmetric with respect to the irreducible linear representation $I_{t}$} of $S$ if $v\in V_{e}^{(I_{t})}$. Similarly, we say that a vector $w\in \mathbb{R}^{m}$ is \emph{symmetric with respect to $I_{t}$} if $w\in V_{i}^{(I_{t})}$.}
\end{defin}

\begin{figure}[ht]
\begin{center}
\begin{tikzpicture}[very thick,scale=1]
\tikzstyle{every node}=[circle, draw=black, fill=white, inner sep=0pt, minimum width=5pt];
        \path (0,0) node (p1) [label = below left: $p_1$] {} ;
        \path (2,0) node (p2) [label = below right: $p_2$] {} ;
        \path (1,2) node (p3) [label = left: $p_3$] {} ;
        \draw (p1)  -- node [draw=white, below left=3pt] {$e_{1}$} (p2);
        \draw (p2) -- node [draw=white, right=2pt] {$e_{3}$} (p3);
        \draw (p3) -- node [draw=white, left=2pt] {$e_{2}$} (p1);
        \draw [dashed, thin] (1,-1) -- (1,3);
        \draw [->, black!60!white, ultra thick] (p1) -- node [draw=white, left=8.5pt] {\setlength{\arraycolsep}{0.1pt}$\left(\begin{array}{r} u_{1}\\u_{2}\end{array}\right)$} (-0.6,0.6);
        \draw [->,  black!60!white, ultra thick] (p2) -- node [draw=white, right=8.5pt] {\setlength{\arraycolsep}{0.1pt}$\left(\begin{array}{r} -u_{1}\\u_{2}\end{array}\right)$} (2.6,0.6);
        \draw [->,  black!60!white, ultra thick] (p3) -- node [draw=white, right=3pt] {\setlength{\arraycolsep}{0.1pt}$\left(\begin{array}{c}0\\u_{3}\end{array}\right)$} (1,2.5);
        \node [draw=white, fill=white] (a) at (1,-1.4) {(a)};
        \end{tikzpicture}
        \hspace{0.02cm}
        \begin{tikzpicture}[very thick,scale=1]
\tikzstyle{every node}=[circle, draw=black, fill=white, inner sep=0pt, minimum width=5pt];
        \path (0,0) node (p1) [label = below left: $p_1$] {} ;
        \path (2,0) node (p2) [label = below right: $p_2$] {} ;
        \path (1,2) node (p3) [label = right: $p_3$] {} ;
        \draw (p1)  -- node [draw=white, below left=3pt] {$e_{1}$} (p2);
        \draw (p2) -- node [draw=white, right=2pt] {$e_{3}$} (p3);
        \draw (p3) -- node [draw=white, left=2pt] {$e_{2}$} (p1);
        \draw [dashed, thin] (1,-1) -- (1,3);
        \draw [->,  black!60!white, ultra thick] (p1) -- node [draw=white, left=8.5pt] {\setlength{\arraycolsep}{0.1pt}$\left(\begin{array}{r} u_{1}\\u_{2}\end{array}\right)$} (-0.6,0.6);
        \draw [->,  black!60!white, ultra thick] (p2)node [draw=white, right=11pt] {\setlength{\arraycolsep}{0.1pt}$\left(\begin{array}{r} u_{1}\\-u_{2}\end{array}\right)$} --  (1.4,-0.6);
        \draw [->,  black!60!white, ultra thick] (p3) -- node [draw=white, left=7pt] {\setlength{\arraycolsep}{0.1pt}$\left(\begin{array}{c}u_{3}\\0\end{array}\right)$} (0.5,2);
        \node [draw=white, fill=white] (b) at (1,-1.4) {(b)};
                \end{tikzpicture}
\end{center}
\vspace{-0.3cm}
\caption{\emph{\emph{(a)} Displacement vectors which are symmetric with respect to $I_1$ (the displacement vector at each joint of $(K_{3},p)$ remains unchanged under $Id$ and $s$); \emph{(b)} displacement vectors which are symmetric with respect to $I_2$ (the displacement vector at each joint of $(K_{3},p)$ remains unchanged under $Id$, but is reversed by $s$).}}
\label{trisymloadres}
\end{figure}
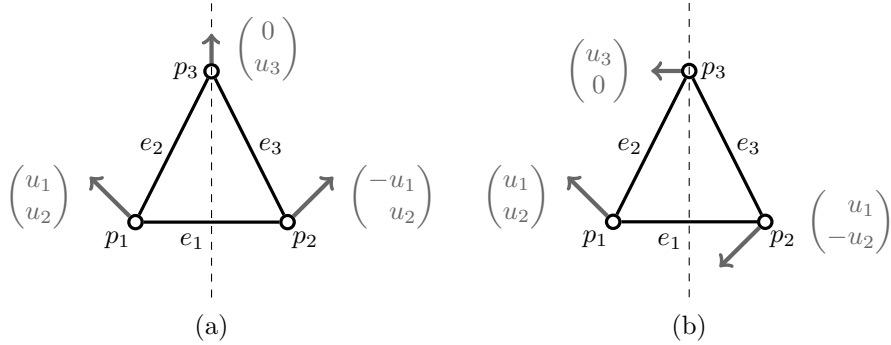

\begin{examp}
\label{blockmatrixexam}
\emph{Let $K_{3}$, $\mathcal{C}_{s}=\{Id,s\}$, and $\Phi$ be as in Example \ref{triangexam}. The symmetry group $\mathcal{C}_{s}$ has two non-equivalent irreducible linear representations, $I_1$ and $I_2$, both of which are of dimension 1 \cite{bishop, Hall}. $I_1$ is the trivial representation which maps both $Id$ and $s$ to the identity transformation, whereas $I_2$ maps $Id$ to the identity transformation and $s$ to the linear transformation $I_2(s)$ which is defined by $I_2(s)(x)=-x$ for all $x\in \mathbb{R}$.
It is easy to see that both of the $H_{e}'$-invariant subspaces $V_{e}^{(I_1)}$ and $V_{e}^{(I_2)}$  of $\mathbb{R}^6$ are of dimension 3. An element of each of these invariant subspaces is shown in Figure \ref{trisymloadres}.
\\\indent Similarly, the $H_{i}'$-invariant subspaces $V_{i}^{(I_1)}$ and $V_{i}^{(I_2)}$  of $\mathbb{R}^3$ are easily found to be of dimension 2 and 1, respectively.
}
\end{examp}

\begin{remark}
\emph{While in Example \ref{blockmatrixexam}, the $H_{e}'$- and $H_{i}'$-invariant subspaces can be found by inspection, this is of course generally not possible. There are, however, some standard methods and algorithms for finding the symmetry adapted bases $B_{e}$ and $B_{i}$ for any given symmetry group. Good sources for these methods are \cite{faessler, weeny}, for example.}
\end{remark}

The following fundamental theorem for analyzing the rigidity properties of a symmetric framework using group representation theory was established in \cite{BS2}.

\begin{theorem}
\label{fundthmrepth} \cite{BS2, BS4}
Let $G$ be a graph, $S$ be a symmetry group with pairwise non-equivalent irreducible linear representations $I_{1},\ldots,I_{r}$, $\Phi$ be a homomorphism from $S$ to $\textrm{Aut}(G)$, and $p\in \bigcap_{x \in S} L_{x,\Phi}$. Then for every $t\in \{1,\ldots,r\}$, we have that if $\mathbf{R}(G,p) u = z$ and $u$ is symmetric with respect to $I_{t}$, then $z$ is also symmetric with respect to $I_{t}$.
\end{theorem}

As an immediate consequence of Theorem \ref{fundthmrepth} we obtain

\begin{cor}
\label{blockdiagon} \cite{KG2, BS2, BS4}
Let $G$ be a graph, $S$ be a symmetry group with pairwise non-equivalent irreducible linear representations $I_{1},\ldots,I_{r}$, $\Phi$ be a homomorphism from $S$ to $\textrm{Aut}(G)$, and $p\in \bigcap_{x \in S} L_{x,\Phi}$. Further, let $T_{e}$ and $T_{i}$ be the matrices of the basis transformations from the canonical bases of the $\mathbb{R}$-vector spaces $\mathbb{R}^{dn}$ and $\mathbb{R}^{m}$ to the bases $B_{e}$ and $B_{i}$, respectively. Then the matrix $\widetilde{\mathbf{R}}(G,p)=T_{i}^{-1}\mathbf{R}(G,p)T_{e}$ is block-diagonalized in such a way that there exists (at most) one submatrix block for each irreducible linear representation $I_{t}$ of $S$.
\end{cor}

Note that in this block-diagonalization of the rigidity matrix $\mathbf{R}(G,p)$, the submatrix block corresponding to $I_{t}$ is a matrix of the size $\textrm{dim }\big(V_{i}^{(I_{t})}\big)\times \textrm{dim }\big(V_{e}^{(I_{t})}\big)$. In particular, a submatrix block can possibly be an `empty matrix' which has rows but no columns or alternatively columns but no rows.\\\indent
Consider the submatrix block $\widetilde{\mathbf{R}}_{t}(G,p)$ of $\widetilde{\mathbf{R}}(G,p)$ corresponding to the irreducible representation $I_t$ of $S$.
By definition, the kernel of  $\widetilde{\mathbf{R}}_{t}(G,p)$ is the space of all infinitesimal motions of $(G,p)$ that are symmetric with respect to $I_t$. Further, $(G,p)$ has a nonzero self-stress that is symmetric with respect to $I_t$ if and only if there exists a linear dependency among the row vectors of $\widetilde{\mathbf{R}}_{t}(G,p)$.
So, the matrix $\widetilde{\mathbf{R}}_{t}(G,p)$ comprises all the information about the infinitesimal motions and self-stresses of $(G,p)$ that are symmetric with respect to $I_t$.\\\indent
Notice that knowledge of only the \emph{size} of the submatrix block $\widetilde{\mathbf{R}}_{t}(G,p)$ already allows us to gain significant insight into the `$I_t$-symmetric'  infinitesimal rigidity properties  of $(G,p)$:
if \begin{equation}\label{maxwell}\textrm{dim }\big(V_{i}^{(I_t)}\big)<\textrm{dim }\big(V_{e}^{(I_t)}\big)-\textrm{dim }\big(W_{e}^{(I_t)}\big)\textrm{,}\end{equation} where $W_{e}^{(I_t)}$ denotes the space of all infinitesimal rigid motions in $V_{e}^{(I_t)}$, then $(G,p)$ clearly possesses an infinitesimal flex that is symmetric with respect to $I_t$. Similarly, if \begin{displaymath}\textrm{dim }\big(V_{i}^{(I_t)}\big)>\textrm{dim }\big(V_{e}^{(I_t)}\big)-\textrm{dim }\big(W_{e}^{(I_t)}\big)\textrm{,}\end{displaymath} then $(G,p)$ clearly possesses a non-zero self-stress which is symmetric with respect to $I_t$ \cite{BS2, BS4}.\\\indent
So, for $(G,p)$ to be isostatic (i.e., infinitesimally rigid and independent), we need to have equality in equation (\ref{maxwell}), for each $t=1 \ldots, r$. These are necessary, but not sufficient conditions for $(G,p)$ to be isostatic, as they cannot predict the presence of paired `equi-symmetric' infinitesimal flexes and self-stresses.\\\indent
The symmetry-extended version of Maxwell's rule described by P. Fowler and S. Guest in \cite{FGsymmax} combines all of these conditions into a single equation using some basic techniques from character theory. With this rule we can detect infinitesimal flexes and self-stresses of $(G,p)$ that are symmetric with respect to  $I_t$ with very little computational effort. In particular, we do not need to determine the bases $B_e$ and $B_i$ in equations (\ref{bases2}) and (\ref{bases2Hi}) explicitly. See also \cite{cfgsw, BS2, BS4}, for further details on this rule.\\\indent
Note that since the dimensions of the subspaces $V_{i}^{(I_t)}$, $V_{e}^{(I_t)}$, and $W_{e}^{(I_t)}$ do not depend on the map $p$ of the framework $(G,p)$, provided that the points $p_1,\ldots, p_n$ span all of $\mathbb{R}^d$, the symmetry-extended version of Maxwell's rule provides information about the infinitesimal rigidity properties for \emph{all  $(S,\Phi)$-generic} realizations of $G$ (see also \cite{BS2,BS4}).

\section{Detection of symmetric flexes}
\label{sec:main}

Our main goal in this section is to find sufficient conditions for the existence of a `symmetry-preserving' flex of a symmetric framework. As we will see, the existence of a flex that preserves some, but not all of the symmetries of a given framework can be predicted in an analogous way.
\\\indent
As usual, we let $G$ be a graph with $V(G)=\{v_{1},\ldots,v_{n}\}$, $S$ be a symmetry group in dimension $d$ with $r$ pairwise non-equivalent irreducible linear representations $I_{1},\ldots,I_{r}$, $\Phi:S\to \textrm{Aut}(G)$ be a homomorphism, and $(G,p)$ be a framework in $\mathscr{R}_{(G,S,\Phi)}$.
In the following, $I_{1}$ will always denote the trivial irreducible linear representation of $S$, i.e., $I_{1}$ denotes the linear representation of dimension one with the property that $I_{1}(x)$ is the identity transformation for all $x\in S$.
\\\indent
Recall from Section \ref{subsec: blockdia} that $V_{e}^{(I_{1})}$ denotes the $H'_{e}$-invariant subspace of $\mathbb{R}^{dn}$ which corresponds to $I_{1}$. So, $p\in V_{e}^{(I_{1})}$ if and only if $H'_{e}(x)(p)=p$ for all $x\in S$. Further, recall from Section \ref{subsec:symgeneric} that if $(G,p)$ is a framework in $\mathscr{R}_{(G,S,\Phi)}$, then $p\in \mathbb{R}^{dn}$ is an element of the subspace $U=\bigcap_{x \in S} L_{x,\Phi}$ of $\mathbb{R}^{dn}$.\\\indent Note that it follows immediately from the definitions of $U$ and $V_{e}^{(I_{1})}$ that $U=V_{e}^{(I_{1})}$, because $p\in U$ if and only if $\mathbf{M^{(x)}}p=\mathbf{P_{\Phi(x)}}p$ for all $x\in S$ if and only if $(\mathbf{P_{\Phi(x)}})^{T}\mathbf{M^{(x)}}p=p$ for all $x\in S$ if and only if $H_{e}(x)p=p$ for all $x\in S$ if and only if $p\in V_{e}^{(I_{1})}$.\\\indent
So, since we are interested in flexes of $(G,p)\in\mathscr{R}_{(G,S,\Phi)}$ that preserve the symmetry of $(G,p)$, we need to restrict the edge functions $f_{G}$ of $G$ and $f_{K_{n}}$ of $K_{n}$ to the subspace $V_{e}^{(I_{1})}$ of $\mathbb{R}^{dn}$. In the following, we let $\tilde{f}_{G}:V_{e}^{(I_{1})}\to \mathbb{R}^{|E(G)|}$ denote the restriction of $f_{G}$ to $V_{e}^{(I_{1})}$, and $\tilde{f}_{K_{n}}:V_{e}^{(I_{1})}\to \mathbb{R}^{\binom{n}{2}}$ denote the restriction of $f_{K_{n}}$ to $V_{e}^{(I_{1})}$. The Jacobian matrices of $\tilde{f}_{G}$ and $\tilde{f}_{K_n}$, evaluated at a point $p\in V_{e}^{(I_{1})}$, are denoted by $d\tilde{f}_{G}(p)$ and $d\tilde{f}_{K_n}(p)$, respectively.

\begin{defin}\label{symregulardef}
\emph{An element $p\in V_{e}^{(I_{1})}$ is said to be a \emph{regular point of $G$ in $V_{e}^{(I_{1})}$} if there exists a neighborhood $N_{p}$ of $p$ in $V_{e}^{(I_{1})}$ so that $\textrm{rank } \big(d\tilde{f}_{G}(p)\big)=\textrm{max }\{\textrm{rank } \big(d\tilde{f}_{G}(q)\big)|\, q\in N_{p}\}$. A \emph{regular point of $K_n$ in $V_{e}^{(I_{1})}$} is defined analogously.}
\end{defin}

\begin{defin}
\label{symmpreservflex}
\emph{An \emph{$(S,\Phi)$-symmetry-preserving flex} of a framework $(G,p)\in \mathscr{R}_{(G,S,\Phi)}$ is a differentiable path $x:[0,1]\to V_{e}^{(I_{1})}$ such that $x(0)=p$ and $x(t)\in \tilde{f}_{G}^{-1}\big(\tilde{f}_{G}(p)\big)\setminus \tilde{f}_{K_{n}}^{-1}\big(\tilde{f}_{K_{n}}(p)\big)$ for all $t\in (0,1]$.}
\end{defin}

\begin{lemma}
\label{hs1}
Let $G$ be a graph, $S$ be a symmetry group, $\Phi:S\to \textrm{Aut}(G)$ be a homomorphism, and $(G,p)$ be a framework in $\mathscr{R}_{(G,S,\Phi)}$.
If $p$ is a regular point of $G$ in $V_{e}^{(I_{1})}$, then there exists a neighborhood $N_{p}$ of $p$ in $V_{e}^{(I_{1})}$ such that $\tilde{f}_{G}^{-1}\big(\tilde{f}_{G}(p)\big)\cap N_{p}$ is a smooth manifold of dimension $\textrm{dim }\big(V_{e}^{(I_{1})}\big)-\textrm{rank }\big(d\tilde{f}_{G}(p)\big)$.
\end{lemma}

\emph{Proof.} The result follows immediately from Proposition 2 (and subsequent remark) in \cite{asiroth}. $\square$

\begin{theorem}
\label{flexthm1}
Let $G$ be a graph with $n$ vertices, $S$ be a symmetry group, $\Phi:S\to \textrm{Aut}(G)$ be a homomorphism, and $(G,p)$ be a framework in $\mathscr{R}_{(G,S,\Phi)}$. If $p$ is a regular point of $G$ in $V_{e}^{(I_{1})}$ and also a regular point of $K_n$ in $V_{e}^{(I_{1})}$, then
\begin{itemize}
\item[(i)] $\textrm{rank }\big(d\tilde{f}_{G}(p)\big)= \textrm{rank }\big(d\tilde{f}_{K_n}(p)\big)$ if and only if $(G,p)$ has no $(S,\Phi)$-symmetry-preserving flex;
\item[(ii)] $\textrm{rank }\big(d\tilde{f}_{G}(p)\big)< \textrm{rank }\big(d\tilde{f}_{K_n}(p)\big)$ if and only if $(G,p)$ has an $(S,\Phi)$-symmetry-preserving flex.
\end{itemize}
\end{theorem}

\emph{Proof.} Since $p$ is a regular point of both $G$ and $K_n$ in $V_{e}^{(I_{1})}$, it follows from Lemma \ref{hs1} that there exist neighborhoods $N_{p}$ and $N'_p$ of $p$ in $V_{e}^{(I_{1})}$ so that $\tilde{f}_{G}^{-1}\big(\tilde{f}_{G}(p)\big)\cap N_{p}$ is a manifold of dimension $\textrm{dim }\big(V_{e}^{(I_{1})}\big)-\textrm{rank }\big(d\tilde{f}_{G}(p)\big)$ and $\tilde{f}_{K_n}^{-1}\big(\tilde{f}_{K_n}(p)\big)\cap N'_{p}$ is a manifold of dimension $\textrm{dim }\big(V_{e}^{(I_{1})}\big)-\textrm{rank }\big(d\tilde{f}_{K_n}(p)\big)$. Since  $\tilde{f}_{K_{n}}^{-1}\big(\tilde{f}_{K_{n}}(p)\big)\cap N''_{p}$ is a submanifold of $\tilde{f}_{G}^{-1}\big(\tilde{f}_{G}(p)\big)\cap N''_{p}$, where $N''_p=N_p\cap N'_p$, it follows that \begin{displaymath}\textrm{rank }\big(d\tilde{f}_{G}(p)\big)\leq \textrm{rank }\big(d\tilde{f}_{K_n}(p)\big)\textrm{.}\end{displaymath} Clearly, $\textrm{rank }\big(d\tilde{f}_{G}(p)\big)= \textrm{rank }\big(d\tilde{f}_{K_n}(p)\big)$ if and only if there exists a neighborhood $N^*_{p}$ of $p$ in $V_{e}^{(I_{1})}$ such that $\tilde{f}_{K_{n}}^{-1}\big(\tilde{f}_{K_{n}}(p)\big)\cap N^*_{p}=\tilde{f}_{G}^{-1}\big(\tilde{f}_{G}(p)\big)\cap N^*_{p}$. Therefore, if $\textrm{rank }\big(d\tilde{f}_{G}(p)\big)= \textrm{rank }\big(d\tilde{f}_{K_n}(p)\big)$, then there does not exist an $(S,\Phi)$-symmetry-preserving flex of $(G,p)$.\\\indent If $\textrm{rank }\big(d\tilde{f}_{G}(p)\big)< \textrm{rank }\big(d\tilde{f}_{K_n}(p)\big)$, then every neighborhood of $p$ in $V_{e}^{(I_{1})}$ contains elements of $\tilde{f}_{G}^{-1}\big(\tilde{f}_{G}(p)\big)\setminus\tilde{f}_{K_{n}}^{-1}\big(\tilde{f}_{K_{n}}(p)\big)$, and hence, by the proof of Proposition 1 in \cite{asiroth} (and references therein), there exists an $(S,\Phi)$-symmetry-preserving flex of $(G,p)$. This completes the proof. $\square$

In order to make further use of Theorem \ref{flexthm1}, we need the following fundamental observations.\\\indent
Recall from Section \ref{subsec: blockdia} that with respect to the bases $B_{e}$ and $B_{i}$, the rigidity matrix of a framework $(G,p)$ in $\mathscr{R}_{(G,S,\Phi)}$ has the block form
\begin{equation}
\label{rigblocks}
\widetilde{\mathbf{R}}(G,p)=\left(\begin{array}{ccc}\widetilde{\mathbf{R}}_{1}(G,p) & & \mathbf{0}\\ & \ddots & \\\mathbf{0} &  & \widetilde{\mathbf{R}}_{r}(G,p) \end{array}\right)\textrm{,}
\end{equation}
where for $t=1,\ldots, r$, the block $\widetilde{\mathbf{R}}_{t}(G,p)$ corresponds to the irreducible linear representation $I_{t}$ of $S$, and the size of the block $\widetilde{\mathbf{R}}_{t}(G,p)$ depends on the dimensions of the subspaces $V_{e}^{(I_{t})}$ of $\mathbb{R}^{dn}$ and $V_{i}^{(I_{t})}$ of $\mathbb{R}^{|E(G)|}$. (In particular, the block $\widetilde{\mathbf{R}}_{t}(G,p)$ is an empty $0\times 0$ matrix if and only if both of the coefficients $\lambda_{t}$ and $\mu_{t}$ in equations (\ref{irrrep}) and (\ref{irrrepHi}) are equal to zero.)\\\indent Since with respect to the bases $B_{e}$ and $B_{i}$, the Jacobian matrix of $f_{G}$, evaluated at $p$, is (up to a constant) the matrix $\widetilde{\mathbf{R}}(G,p)$, it follows that with respect to the bases $B_{e}^{(I_{1})}$ and $B_{i}$, the Jacobian matrix of $\tilde{f}_{G}$, evaluated at the point $p\in V_{e}^{(I_{1})}$, is (up to a constant) the matrix
\begin{displaymath}
\left(\begin{array}{c} \widetilde{\mathbf{R}}_{1}(G,p)\\ \mathbf{0}\\ \vdots\\ \mathbf{0}\end{array}\right)\textrm{.}
\end{displaymath}
Thus, we have \begin{equation}\label{ranks1}\textrm{rank }\big(\widetilde{\mathbf{R}}_{1}(G,p)\big)=\textrm{rank }\big(d\tilde{f}_{G}(p)\big)\textrm{.}\end{equation}
\indent Furthermore, note that if $K_{n}$ is the complete graph on the vertex set $V(G)$, then with respect to the bases $B_{e}$ and $\widehat{B}_{i}$, where $\widehat{B}_{i}$ is an appropriate extension of the basis $B_{i}$, the rigidity matrix of $(K_{n},p)$ has a block form analogous to the one of $\widetilde{\mathbf{R}}(G,p)$ in (\ref{rigblocks}), namely
\begin{equation}
\label{rigblocks2}
\widetilde{\mathbf{R}}(K_{n},p)=\left(\begin{array}{ccc}\widetilde{\mathbf{R}}_{1}(K_{n},p) & & \mathbf{0}\\ & \ddots & \\\mathbf{0} &  & \widetilde{\mathbf{R}}_{r}(K_{n},p) \end{array}\right)\textrm{.}\nonumber
\end{equation}
Clearly, $\widetilde{\mathbf{R}}_{t}(G,p)$ is a submatrix of $\widetilde{\mathbf{R}}_{t}(K_{n},p)$ for all $t=1,\ldots, r$. Moreover, analogously to (\ref{ranks1}), we have \begin{equation}\label{ranks2}\textrm{rank }\big(\widetilde{\mathbf{R}}_{1}(K_n,p)\big)=\textrm{rank }\big(d\tilde{f}_{K_n}(p)\big)\textrm{.}\end{equation}

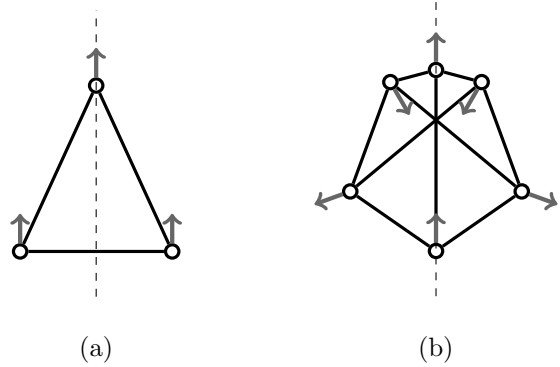
\begin{figure}[htp]
\begin{center}
\begin{tikzpicture}[very thick,scale=1]
\tikzstyle{every node}=[circle, draw=black, fill=white, inner sep=0pt, minimum width=5pt];
        \node  (p1) at (0,0) {} ;
        \node  (p2) at (2,0) {} ;
       \node  (p3) at (1,2.2) {} ;
        \draw (p1)  --  (p2);
        \draw (p3)  --  (p2);
        \draw (p1)  --  (p3);
        \draw [dashed, thin] (1,-0.6) -- (1,3.2);
        \draw [->, black!60!white, ultra thick] (p1) --  (0,0.5);
        \draw [->,black!60!white, ultra thick] (p2) --  (2,0.5);
        \draw [->, black!60!white, ultra thick] (p3) --  (1,2.7);
        \node [draw=white, fill=white] (a) at (1,-1.3) {(a)};
        \end{tikzpicture}
        \hspace{1.5cm}
        \begin{tikzpicture}[very thick,scale=1]
\tikzstyle{every node}=[circle, draw=black, fill=white, inner sep=0pt, minimum width=5pt];
        \node  (p1) at (60:1.2cm) {} ;
        \node  (p2) at (90:1.2cm) {} ;
        \node  (p3) at (120:1.2cm) {} ;
        \node  (p4) at (200:1.2cm) {} ;
        \node  (p5) at (270:1.2cm) {} ;
        \node  (p6) at (340:1.2cm) {} ;
        \draw (p1)  --  (p2);
        \draw (p3)  --  (p2);
        \draw (p3)  --  (p4);
        \draw (p4)  --  (p5);
        \draw (p5)  --  (p6);
        \draw (p1)  --  (p6);
        \draw (p5)  --  (p2);
        \draw (p3)  --  (p6);
        \draw (p1)  --  (p4);
        \draw [dashed, thin] (270:1.8cm) -- (90:2.1cm);
        \draw [->, black!60!white, ultra thick] (p1) --  (60:0.7cm);
        \draw [->, black!60!white, ultra thick] (p3) --  (120:0.7cm);
        \draw [->, black!60!white, ultra thick] (p2) --  (90:1.7cm);
        \draw [->, black!60!white, ultra thick] (p4) --  (200:1.7cm);
        \draw [->, black!60!white, ultra thick] (p5) --  (270:0.7cm);
        \draw [->, black!60!white, ultra thick] (p6) --  (340:1.7cm);
        \node [draw=white, fill=white] (b) at (270:2.5cm) {(b)};
        \end{tikzpicture}
         \end{center}
        \vspace{-0.3cm}
        \caption{\emph{Fully $(S,\Phi)$-symmetric infinitesimal motions of frameworks: \emph{(a)} a fully $(\mathcal{C}_{s},\Phi)$-symmetric infinitesimal rigid motion of $(K_{3},p)\in \mathscr{R}_{(K_{3},\mathcal{C}_{s},\Phi)}$; \emph{(b)} a fully $(\mathcal{C}_{s},\Phi)$-symmetric infinitesimal flex of $(K_{3,3},p)\in \mathscr{R}_{(K_{3,3},\mathcal{C}_{s},\Phi)}$. Since each of the above frameworks is an injective realization, the type $\Phi$ is uniquely determined in each case \cite{BS4,BS1}.}}
\label{fulsym}
        \end{figure}

Recall from Definition \ref{usym}
that if $u\in V_{e}^{(I_{1})}$, then $u$ is said to be symmetric with respect to $I_{1}$.
So, if we think of the vector $u\in \mathbb{R}^{dn}$ as a set of displacement vectors with one vector at each joint of $(G,p)\in \mathscr{R}_{(G,S,\Phi)}$, then $u$ is symmetric with respect to $I_{1}$ if and only if all of the displacement vectors remain unchanged under all symmetry operations in $S$.
A vector $u\in \mathbb{R}^{dn}$ that is symmetric with respect to $I_{1}$ can therefore also be termed \emph{fully $(S,\Phi)$-symmetric} \cite{FG4, KG1} (see also Figure \ref{fulsym}).

\begin{theorem}
\label{flexthm2}
Let $G$ be a graph, $S$ be a symmetry group in dimension $d$, $\Phi:S\to \textrm{Aut}(G)$ be a homomorphism, and $(G,p)$ be a framework in $\mathscr{R}_{(G,S,\Phi)}$ with the property that the points $p(v)$, $v\in V(G)$, span all of $\mathbb{R}^{d}$. If $p$ is a regular point of $G$ in $V_{e}^{(I_{1})}$ and also a regular point of $K_n$ in $V_{e}^{(I_{1})}$ and there exists a fully $(S,\Phi)$-symmetric infinitesimal flex of $(G,p)$, then there also exists an $(S,\Phi)$-symmetry-preserving flex of $(G,p)$.
\end{theorem}

\emph{Proof.} Let $K_{n}$ be the complete graph on $V(G)$. Since the points $p(v)$, $v\in V(G)$, span all of $\mathbb{R}^{d}$, the kernel of $\mathbf{R}(K_n,p)$ is the space of all infinitesimal rigid motions of $(G,p)$ and the kernel of $\widetilde{\mathbf{R}}_{1}(K_{n},p)$ is the space of all fully $(S,\Phi)$-symmetric infinitesimal rigid motions of $(G,p)$. So, since  $(G,p)$ has an infinitesimal flex which is fully $(S,\Phi)$-symmetric, we have $\textrm{nullity }\big(\widetilde{\mathbf{R}}_{1}(G,p)\big)>\textrm{nullity }\big(\widetilde{\mathbf{R}}_{1}(K_{n},p)\big)$, and hence \begin{equation}\label{rankinproof}\textrm{rank }\big(\widetilde{\mathbf{R}}_{1}(G,p)\big)<\textrm{rank }\big(\widetilde{\mathbf{R}}_{1}(K_{n},p)\big)\textrm{.}\end{equation}
Since, by (\ref{ranks1}), we have $\textrm{rank }\big(\widetilde{\mathbf{R}}_{1}(G,p)\big)=\textrm{rank }\big(d\tilde{f}_{G}(p)\big)$ and, by (\ref{ranks2}), we have $\textrm{rank }\big(\widetilde{\mathbf{R}}_{1}(K_n,p)\big)=\textrm{rank }\big(d\tilde{f}_{K_n}(p)\big)$, it follows from (\ref{rankinproof}) that \begin{displaymath}\textrm{rank }\big(d\tilde{f}_{G}(p)\big)< \textrm{rank }\big(d\tilde{f}_{K_n}(p)\big)\textrm{.}\end{displaymath} The result now follows from Theorem \ref{flexthm1}. $\square$

The above results concerning the subspace $V_{e}^{(I_{1})}$ of $\mathbb{R}^{dn}$ may be transferred analogously
to the affine subspaces of $\mathbb{R}^{dn}$ of the form $p+V_{e}^{(I_{t})}$, where $t\neq 1$.
\\\indent More precisely, if we define a point $q\in p+V_{e}^{(I_{t})}$ to be a \emph{regular point of a graph $G$ in $p+V_{e}^{(I_{t})}$} if there exists a neighborhood $N_{q}$ of $q$ in $p+V_{e}^{(I_{t})}$ so that $\textrm{rank } \big(d\tilde{f}_{G}(q)\big)=\textrm{max }\{\textrm{rank } \big(d\tilde{f}_{G}(q')\big)|\, q'\in N_{q}\}$, where $\tilde{f}_{G}$ denotes the restriction of the edge function $f_{G}$ to $p+V_{e}^{(I_{t})}$, then the following results can be proved completely analogously to the Theorems \ref{flexthm1} and \ref{flexthm2}.

\begin{theorem}
\label{flexthm1b}
Let $G$ be a graph with $n$ vertices, $S$ be a symmetry group, $\Phi:S\to \textrm{Aut}(G)$ be a homomorphism, and $(G,p)$ be a framework in $\mathscr{R}_{(G,S,\Phi)}$. If $p$ is a regular point of $G$ in $p+V_{e}^{(I_{t})}$ and also a regular point of $K_n$ in $p+V_{e}^{(I_{t})}$, then
\begin{itemize}
\item[(i)] $\textrm{rank }\big(d\tilde{f}_{G}(p)\big)= \textrm{rank }\big(d\tilde{f}_{K_n}(p)\big)$ if and only if $(G,p)$ does not have a flex $x$ with $x(t)\in p+V_{e}^{(I_{t})}$ for all $t\in[0,1]$;
\item[(ii)] $\textrm{rank }\big(d\tilde{f}_{G}(p)\big)< \textrm{rank }\big(d\tilde{f}_{K_n}(p)\big)$ if and only if $(G,p)$ has a flex $x$ with $x(t)\in p+V_{e}^{(I_{t})}$ for all $t\in[0,1]$.
\end{itemize}
\end{theorem}

\begin{theorem}
\label{flexthm2b}
Let $G$ be a graph, $S$ be a symmetry group in dimension $d$, $\Phi:S\to \textrm{Aut}(G)$ be a homomorphism, and $(G,p)$ be a framework in $\mathscr{R}_{(G,S,\Phi)}$ with the property that the points $p(v)$, $v\in V(G)$, span all of $\mathbb{R}^{d}$. If $p$ is a regular point of $G$ in $p+V_{e}^{(I_{t})}$ and also a regular point of $K_n$ in $p+V_{e}^{(I_{t})}$ and there exists an infinitesimal flex $u$ of $(G,p)$ with $u\in V_{e}^{(I_{t})}$, then there also exists a flex $x$ of $(G,p)$ with $x(t)\in p+V_{e}^{(I_{t})}$ for all $t\in[0,1]$.
\end{theorem}

Note that if we define $\textrm{ker }(I_{t})=\{x\in S|\, I_{t}(x)=id\}$, where $id$ is the identity transformation, then $\textrm{ker }(I_{t})$ is a normal subgroup of $S$ (see \cite{liebeck}, for example). Therefore, Theorems \ref{flexthm1b} and \ref{flexthm2b} provide us with sufficient conditions for the existence of a flex of $(G,p)$ that preserves the sub-symmetry of $(G,p)$ given by $\textrm{ker }(I_{t})$ and $\Phi|_{\textrm{ker }(I_{t})}$.

An important property of the subspace $V_{e}^{(I_{1})}$ which does not hold for the affine subspaces $p+V_{e}^{(I_{t})}$, where $t\neq 1$, is that, by Corollary \ref{blockdiagon}, for every $q\in V_{e}^{(I_{1})}$, the rigidity matrix $\widetilde{\mathbf{R}}(G,q)$ has the same block structure as the rigidity matrix $\widetilde{\mathbf{R}}(G,p)$. Thus, $p\in V_{e}^{(I_{1})}$ is a regular point of $G$ in $V_{e}^{(I_{1})}$ if and only if there exists a neighborhood $N_{p}$ of $p$ in $V_{e}^{(I_{1})}$ so that $\textrm{rank }\big(\widetilde{\mathbf{R}}_{1}(G,p)\big)\geq \textrm{rank }\big(\widetilde{\mathbf{R}}_{1}(G,q)\big)$ for all $q\in N_{p}$.\\\indent Similarly, $p\in V_{e}^{(I_{1})}$ is a regular point of $K_n$ in $V_{e}^{(I_{1})}$ if and only if there exists a neighborhood $N_{p}$ of $p$ in $V_{e}^{(I_{1})}$ so that $\textrm{rank }\big(\widetilde{\mathbf{R}}_{1}(K_n,p)\big)\geq \textrm{rank }\big(\widetilde{\mathbf{R}}_{1}(K_n,q)\big)$ for all $q\in N_{p}$.\\\indent The fact that regular points of $G$ and $K_n$ in $V_{e}^{(I_{1})}$ can be characterized in this way is essential to proving all the remaining results of this section. These results turn out to be very useful for  practical applications of Theorem \ref{flexthm2}, as we will see in Section \ref{sec:ex} (see also \cite{SchW} as well as Section 6.3 in \cite{BS4}).

\begin{theorem}
\label{spanandreg}
Let $G$ be a graph with $n$ vertices, $S$ be a symmetry group in dimension $d$, $\Phi:S\to \textrm{Aut}(G)$ be a homomorphism, and $(G,p)$ be a framework in $\mathscr{R}_{(G,S,\Phi)}$. If the points $p(v)$, $v\in V(G)$, span all of $\mathbb{R}^{d}$, then $p$ is a regular point of $K_n$ in $V_{e}^{(I_{1})}$.
\end{theorem}

\emph{Proof.} Since the points $p(v)$, $v\in V(G)$, span all of $\mathbb{R}^{d}$, there exists a neighborhood $N_p$ of $p$ in $V_{e}^{(I_{1})}$ so that for all $q\in N_p$, the points $q(v)$, $v\in V(G)$, also span all of $\mathbb{R}^{d}$. Therefore, for all $q\in N_p$, the dimension of the subspace of $\mathbb{R}^{dn}$ consisting of all fully $(S,\Phi)$-symmetric infinitesimal rigid motions of $(G,p)$ is equal to the dimension of the subspace of $\mathbb{R}^{dn}$ consisting of all fully $(S,\Phi)$-symmetric infinitesimal rigid motions of $(G,q)$ (see \cite{BS2, BS4} for details). Therefore, we have $\textrm{rank }\big(\widetilde{\mathbf{R}}_{1}(K_n,p)\big)=\textrm{rank }\big(\widetilde{\mathbf{R}}_{1}(K_{n},q)\big)$ or equivalently, by (\ref{ranks2}), $\textrm{rank }\big(d\tilde{f}_{K_n}(p)\big)= \textrm{rank }\big(d\tilde{f}_{K_n}(q)\big)$ for all $q\in N_p$. Thus, $p$ is a regular point of $K_n$ in $V_{e}^{(I_{1})}$.  $\square$

By Theorem \ref{spanandreg}, the condition that $p$ is a regular point of $K_n$ in $V_{e}^{(I_{1})}$ may be omitted in Theorem \ref{flexthm2}.

\begin{theorem}
\label{reggen}
Let $G$ be a graph, $S$ be a symmetry group, $\Phi:S\to \textrm{Aut}(G)$ be a homomorphism, and $(G,p)$ be a framework in $\mathscr{R}_{(G,S,\Phi)}$. If $p$ is $(S,\Phi)$-generic, then $p$ is a regular point of $G$ in $V_{e}^{(I_{1})}$.
\end{theorem}

\emph{Proof.} Suppose $G$ is a graph with $n$ vertices and $S$ is a symmetry group in dimension $d$ with $r$ pairwise non-equivalent irreducible representations $I_{1},\ldots,I_{r}$. Fix a basis $\mathscr{B}_{U}=\{u_{1},\ldots,u_{k}\}$ of $U=V_{e}^{(I_{1})}=\bigcap_{x \in S}L_{x,\Phi}$ and let $p=t_{1}u_{1}+\ldots+t_{k}u_{k}$. Then the symmetry-adapted `indeterminate' rigidity matrix $\mathbf{R}_{\mathscr{B}_{U}}(n,d)$ for $\mathscr{R}_{(G,S,\Phi)}$  is a matrix in the variables $t'_{1},\ldots,t'_{k}$. More precisely, the entries of $\mathbf{R}_{\mathscr{B}_{U}}(n,d)$ are elements of the quotient field of the integral domain $\mathbb{R}[ t'_{1},\ldots,t'_{k}]$. Over this field we can again do linear algebra. We let $\mathbf{R}^{(G)}_{\mathscr{B}_{U}}(n,d)$ denote the submatrix of $\mathbf{R}_{\mathscr{B}_{U}}(n,d)$ that corresponds to the submatrix $\mathbf{R}(G,p)$ of $\mathbf{R}(K_{n},p)$, i.e., $\mathbf{R}^{(G)}_{\mathscr{B}_{U}}(n,d)$ is obtained from $\mathbf{R}_{\mathscr{B}_{U}}(n,d)$ by deleting those rows that do not correspond to edges of $G$. If we replace each variable $t'_{i}$ in $\mathbf{R}^{(G)}_{\mathscr{B}_{U}}(n,d)$ with $t_{i}$, then, by Remark \ref{symindetrigmatrem}, we obtain the rigidity matrix $\mathbf{R}(G,p)$. Therefore, \begin{displaymath}\textrm{rank }\big(\mathbf{R}(G,p)\big)\leq \textrm{rank }\big(\mathbf{R}^{(G)}_{\mathscr{B}_{U}}(n,d)\big)\textrm{.}\end{displaymath} Since $(G,p)$ is $(S,\Phi)$-generic, we also have \begin{displaymath}\textrm{rank }\big(\mathbf{R}(G,p)\big)\geq \textrm{rank }\big(\mathbf{R}^{(G)}_{\mathscr{B}_{U}}(n,d)\big)\textrm{,}\end{displaymath} and hence \begin{equation}\label{ranks}\textrm{rank }\big(\mathbf{R}(G,p)\big)= \textrm{rank }\big(\mathbf{R}^{(G)}_{\mathscr{B}_{U}}(n,d)\big)\textrm{.}\end{equation}
Now, let $T_{e}$ be the matrix of the basis transformation from the canonical basis of the $\mathbb{R}$-vector space $\mathbb{R}^{dn}$ to the basis $B_{e}$, and let $T_{i}$ be the matrix of the basis transformation from the canonical basis of the $\mathbb{R}$-vector space $\mathbb{R}^{|E(G)|}$ to the basis $B_{i}$, so that the matrix $\widetilde{\mathbf{R}}(G,p)=T_{i}^{-1}\mathbf{R}(G,p)T_{e}$ is block-diagonalized as in (\ref{rigblocks}). Then, by Corollary \ref{blockdiagon}, the matrix $\widetilde{\mathbf{R}}^{(G)}_{\mathscr{B}_{U}}(n,d)=T_{i}^{-1}\mathbf{R}^{(G)}_{\mathscr{B}_{U}}(n,d)T_{e}$ has the same block form as $\widetilde{\mathbf{R}}(G,p)$. For $t=1,\ldots,r$, let $\widetilde{\mathbf{R}}^{(G)}_{t}(n,d)$ denote the block of $\widetilde{\mathbf{R}}^{(G)}_{\mathscr{B}_{U}}(n,d)$ that corresponds to the block $\widetilde{\mathbf{R}}_{t}(G,p)$ of $\widetilde{\mathbf{R}}(G,p)$. Since the rank of a matrix is invariant under a basis transformation, and since the rank of a block-diagonalized matrix is equal to the sum of the ranks of its blocks, it follows from equation (\ref{ranks}) that
\begin{eqnarray} \sum_{t=1}^{r}\textrm{rank }\big(\widetilde{\mathbf{R}}_{t}(G,p)\big)&=&\textrm{rank }\big(\widetilde{\mathbf{R}}(G,p)\big)\nonumber\\
&=&\textrm{rank }\big(\mathbf{R}(G,p)\big)\nonumber\\
&=&\textrm{rank }\big(\mathbf{R}^{(G)}_{\mathscr{B}_{U}}(n,d)\big)\nonumber\\
&=&\textrm{rank }\big(\widetilde{\mathbf{R}}^{(G)}_{\mathscr{B}_{U}}(n,d)\big)=\sum_{t=1}^{r}\textrm{rank }\big(\widetilde{\mathbf{R}}^{(G)}_{t}(n,d)\big)\textrm{.}\nonumber\end{eqnarray}
Since we clearly have $\textrm{rank }\big(\widetilde{\mathbf{R}}_{t}(G,p)\big)\leq\textrm{rank }\big(\widetilde{\mathbf{R}}^{(G)}_{t}(n,d)\big)$ for each $t$, it follows that $\textrm{rank }\big(\widetilde{\mathbf{R}}_{t}(G,p)\big)=\textrm{rank }\big(\widetilde{\mathbf{R}}^{(G)}_{t}(n,d)\big)$ for each $t$. This gives the result. $\square$

\begin{cor} \label{genflexcor} Let $G$ be a graph, $S$ be a symmetry group in dimension $d$, $\Phi:S\to \textrm{Aut}(G)$ be a homomorphism, and $(G,p)$ be a framework in $\mathscr{R}_{(G,S,\Phi)}$ with the property that the points $p(v)$, $v\in V(G)$, span all of $\mathbb{R}^{d}$. If $(G,p)$ is $(S,\Phi)$-generic and $(G,p)$ has a fully $(S,\Phi)$-symmetric infinitesimal flex, then there also exists an $(S,\Phi)$-symmetry-preserving flex of $(G,p)$.
\end{cor}

\emph{Proof.} The result follows immediately from Theorems \ref{flexthm2}, \ref{spanandreg}, and \ref{reggen}. $\square$

In Section \ref{sec:ex}, we will use Corollary \ref{genflexcor} to prove the existence of an $(S,\Phi)$-symmetry-preserving flex for a variety of symmetric octahedral frameworks in 3-space (including two of the three types of `Bricard octahedra' \cite{bricard}). A number of other classes of frameworks, such as the ones examined in \cite{bottema, consusp, FG4, BS4, tarnaiely}, can also be shown to be fexible with the help of Corollary \ref{genflexcor}. \\\indent Note that the framework in Figure \ref{fulsym} (b) is not $(\mathcal{C}_s,\Phi)$-generic (recall Section \ref{subsec:symgeneric} and Figure \ref{K33types}), so that Corollary \ref{genflexcor} does not apply to this framework.
In fact, it can be verified that the framework in Figure \ref{fulsym} (b) does not possess any flex, let alone a $(\mathcal{C}_s,\Phi)$-symmetry-preserving flex.
\\\indent
Corollary \ref{genflexcor} is a symmetrized version of Corollary \ref{asigen} in Section \ref{subsec:basicgenres}. Next, we show that a symmetrized version of Corollary \ref{asiind} can also be obtained from the previous results.

\begin{cor}\label{indepfullysymfini} Let $G$ be a graph, $S$ be a symmetry group in dimension $d$, $\Phi:S\to \textrm{Aut}(G)$ be a homomorphism, and $(G,p)$ be a framework in $\mathscr{R}_{(G,S,\Phi)}$ with the property that the points $p(v)$, $v\in V(G)$, span all of $\mathbb{R}^{d}$. If the block $\widetilde{\mathbf{R}}_{1}(G,p)$ of the block-diagonalized rigidity matrix $\widetilde{\mathbf{R}}(G,p)$ has linearly independent rows and $(G,p)$ has a fully $(S,\Phi)$-symmetric infinitesimal flex, then there also exists an $(S,\Phi)$-symmetry-preserving flex of $(G,p)$.
\end{cor}

\emph{Proof.} Since the block matrix $\widetilde{\mathbf{R}}_{1}(G,p)$ has linearly independent rows, $p$ is a regular point of $G$ in $V_{e}^{(I_{1})}$. The result now follows from Theorems \ref{flexthm2} and \ref{spanandreg}. $\square$

Corollary \ref{indepfullysymfini} confirms the observation made by R. Kangwai and S. Guest in \cite{KG1}. Note that the condition that the block matrix $\widetilde{\mathbf{R}}_{1}(G,p)$ has linearly independent rows is equivalent to the condition that the framework $(G,p)$ has no fully $(S,\Phi)$-symmetric non-zero self-stress, i.e., a non-zero self-stress in the subspace $V_{i}^{(I_{1})}$ of $\mathbb{R}^{|E(G)|}$. In particular, it follows that if $(G,p)$ is independent (i.e., $(G,p)$ does not possess \emph{any} non-zero self-stress) and there exists a fully $(S,\Phi)$-symmetric infinitesimal flex of $(G,p)$, then there also exists an $(S,\Phi)$-symmetry-preserving flex of $(G,p)$.
\\\indent
In order to apply Corollary \ref{indepfullysymfini} to a given framework $(G,p)$, we need to compute the rank of the submatrix block $\widetilde{\mathbf{R}}_{1}(G,p)$. This can be done by finding the block-diagonalized rigidity matrix $\widetilde{\mathbf{R}}(G,p)$ with the methods and algorithms described in \cite{faessler, weeny}, for example.
\\\indent
The rank of the submatrix block $\widetilde{\mathbf{R}}_{1}(G,p)$ can also be determined directly by finding the rank of an appropriate `orbit rigidity matrix' whose columns and rows correspond to a set of representatives for the orbits of the group action from $S\times V(G)$ to $V(G)$ that sends $(x,v)$ to $\Phi(x)(v)$ and a set of  representatives for the orbits of the group action from $S\times E(G)$ to $E(G)$ that sends $(x,e)$ to $\Phi(x)(e)$, respectively. The kernel of this matrix is the space of fully $(S,\Phi)$-symmetric infinitesimal motions of $(G,p)$ and the cokernel of this matrix is the  space of fully $(S,\Phi)$-symmetric self-stresses of $(G,p)$. Further details on the `orbit rigidity matrix' will be presented in \cite{SchW}.

\section{Examples}
\label{sec:ex}

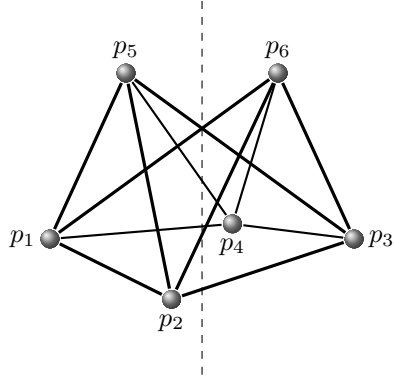
\begin{figure}[htp]
\begin{center}
\begin{tikzpicture}[very thick,scale=1]
\tikzstyle{every node}=[circle, fill=white, inner sep=0pt, minimum width=5pt];
\draw [dashed, thin] (0,-1.5) -- (0,3.5);
\linespread{1.0}
\node [circle, shade, ball color=black!40!white, inner sep=0pt, minimum width=7pt](p1) at (-2,0.3) {};
\node [circle, shade, ball color=black!40!white, inner sep=0pt, minimum width=7pt](p3) at (2,0.3) {};
\node [circle, shade, ball color=black!40!white, inner sep=0pt, minimum width=7pt](p2) at (-0.4,-0.5) {};
\node [circle, shade, ball color=black!40!white, inner sep=0pt, minimum width=7pt](p4) at (0.4,0.5) {};
\node [circle, shade, ball color=black!40!white, inner sep=0pt, minimum width=7pt](p5) at (-1,2.5) {};
\node [circle, shade, ball color=black!40!white, inner sep=0pt, minimum width=7pt](p6) at (1,2.5) {};
\draw (p1) -- (p2)node[rectangle, draw=white, anchor=south, below=5pt] {$p_2$};
\draw (p2) -- (p3) node[rectangle, draw=white,anchor=south, right=5pt] {$p_3$};
\draw [thick](p3) -- (p4) node[rectangle, draw=white,anchor=south, below=5pt] {$p_4$};
\draw [thick](p4) -- (p1) node[rectangle, draw=white,anchor=south, left=5pt] {$p_1$};
\draw (p5)node[rectangle, draw=white,anchor=south, above=5pt] {$p_5$} -- (p1);
\draw (p5) -- (p2);
\draw (p5) -- (p3);
\draw [thick](p5) -- (p4);
\draw (p6)node[rectangle, draw=white,anchor=south, above=5pt] {$p_6$} -- (p1);
\draw (p6) -- (p2);
\draw (p6) -- (p3);
\draw [thick](p6) -- (p4);
\end{tikzpicture}
\end{center}
\vspace{-0.3cm}
\caption{\emph{Flexible Bricard octahedron with point group $\mathcal{C}_2$.}}
\label{figurebricardc2}
\end{figure}

In his famous paper from 1897, the French engineer R. Bricard proved that if an octahedron in 3-space with no self-intersecting faces is realized as a framework by placing bars along edges, and joints at vertices, then this framework must be rigid \cite{bricard}. Moreover, he showed that there exist three distinct types
of octahedra with self-intersecting faces whose realizations as frameworks are flexible.
Two of these three types of octahedra possess non-trivial symmetries: Bricard octahedra of the first type have a half-turn symmetry and Bricard octahedra of the second type have a mirror symmetry.
In the following, we consider both of these types of symmetric Bricard octahedra (as well as octahedra with dihedral symmetry) and use the results of Section \ref{sec:main} to not only show that they are flexible, but also that they possess a `symmetry-preserving' flex.\\\indent Various other treatments of the Bricard octahedra can be found in \cite{Baker,Stachel}, for example. R. Connelly's celebrated counterexample to Euler's rigidity conjecture from 1776 (see \cite{euler})
is also based on a flexible Bricard octahedron (of the first type) \cite{concounter, con1}.
\\\indent
Let $G$ be the graph of the octahedron (see Figure \ref{figurebricardc2}), $\mathcal{C}_{2}$ be a `half-turn' symmetry group in dimension 3, and $\Phi_a:\mathcal{C}_{2}\to \textrm{Aut}(G)$ be the homomorphism defined by
\begin{eqnarray}\Phi_a(Id)& = &id\nonumber\\ \Phi_a(C_2) &=& (v_{1}\,v_{3})(v_{2}\,v_{4})(v_{5}\,v_{6})\textrm{.}\nonumber
\end{eqnarray}
From the symmetry-extended version of Maxwell's rule, applied to the $(\mathcal{C}_{2},\Phi_a)$-generic framework $(G,p)$ in Figure \ref{figurebricardc2}, $\mathcal{C}_{2}$, and $\Phi_a$, we may deduce that  $(G,p)$ has a fully $(\mathcal{C}_{2},\Phi_a)$-symmetric infinitesimal flex \cite{FGsymmax, BS2, BS4}. It follows from Corollary \ref{genflexcor} that $(G,p)$ also has a $(\mathcal{C}_{2},\Phi_a)$-symmetry-preserving flex.
\\\indent Since the symmetry-extended version of Maxwell's rule yields the same counts for \emph{any} $(\mathcal{C}_{2},\Phi_a)$-generic realization of $G$, any such realization has a $(\mathcal{C}_{2},\Phi_a)$-symmetry-preserving flex.

\begin{remark}
\emph{Note that some of the configurations that lie on the path of the $(\mathcal{C}_{2},\Phi_a)$-symmetry-preserving flex of $(G,p)$ are not $(\mathcal{C}_{2},\Phi_a)$-generic.\\\indent For example, the $(\mathcal{C}_{2},\Phi_a)$-symmetry-preserving flex of $(G,p)$ passes through a configuration $q$ with the property that the four points $q_1, q_2, q_3$, and $ q_4$ are coplanar. The framework $(G,q)$ is therefore clearly not $(\mathcal{C}_{2},\Phi_a)$-generic. However, by computing the rank of $\widetilde{\mathbf{R}}_{1}(G,q)$ and showing that it is equal to the rank of $\widetilde{\mathbf{R}}_{1}(G,p)$, where $p$ is $(\mathcal{C}_{2},\Phi_a)$-generic, the configuration $q$ can be proven to be a regular point of $G$ in $V_{e}^{(I_{1})}$, where $I_1$ is the trivial irreducible representation of $\mathcal{C}_{2}$. So, Theorem \ref{flexthm2} can be used in this case to prove the existence of a $(\mathcal{C}_{2},\Phi_a)$-symmetry-preserving flex of $(G,q)$.
}
\end{remark}

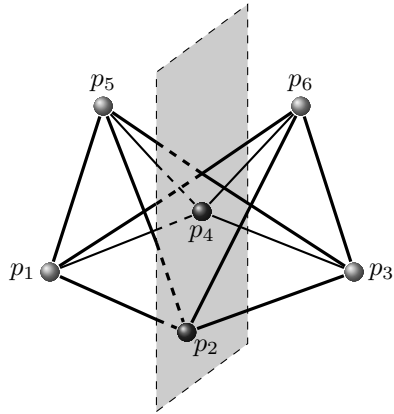
\begin{figure}[htp]
\begin{center}
\begin{tikzpicture}[very thick,scale=1]
\tikzstyle{every node}=[circle, fill=white, inner sep=0pt, minimum width=5pt];
\filldraw[fill=black!20!white, draw=black, thin, dashed]
    (-0.6,-1.55) -- (0.6,-0.65) -- (0.6,3.85) -- (-0.6,2.95) -- (-0.6,-1.55);
\node [circle, shade, ball color=black!40!white, inner sep=0pt, minimum width=7pt](p1) at (-2,0.3) {};
\node [circle, shade, ball color=black!40!white, inner sep=0pt, minimum width=7pt](p3) at (2,0.3) {};
\node [circle, shade, ball color=black!80!white, inner sep=0pt, minimum width=7pt](p2) at (-0.2,-0.5) {};
\node [circle, shade, ball color=black!80!white, inner sep=0pt, minimum width=7pt](p4) at (0,1.1) {};
\node [circle, shade, ball color=black!40!white, inner sep=0pt, minimum width=7pt](p5) at (-1.3,2.5) {};
\node [circle, shade, ball color=black!40!white, inner sep=0pt, minimum width=7pt](p6) at (1.3,2.5) {};
\draw (p1) -- (-0.6,-0.33);
\draw [dashed] (-0.6,-0.33) -- (p2)node[rectangle, fill=black!20!white,anchor=south, below right=2pt] {$p_2$};
\draw [thick](p3) -- (p4) node[rectangle, fill=black!20!white, anchor=south, below=5pt] {$p_4$};
\draw (p2) -- (p3) node[rectangle, draw=white,anchor=south, right=5pt] {$p_3$};
\draw [dashed, thick](p4) -- (-0.6,0.86);
\draw [thick](-0.6,0.86) -- (p1) node[rectangle, draw=white,anchor=south, left=5pt] {$p_1$};
\draw (p5)node[rectangle, draw=white,anchor=south, above=5pt] {$p_5$} -- (p1);
\draw (p5) -- (-0.6,0.69);
\draw [dashed](-0.6,0.69) -- (p2);
\draw (p5) -- (-0.6,2.02);
\draw [dashed](-0.6,2.02)--(0,1.6053);
\draw (0,1.6053) -- (p3);
\draw [thick](p5) -- (-0.6,1.73);
\draw [dashed, thick] (-0.6,1.73) -- (p4);
\draw (p6) node[rectangle, draw=white,anchor=south, above=5pt] {$p_6$}-- (0,1.605);
\draw [dashed](0,1.605)--(-0.6,1.205);
\draw (-0.605,1.205) -- (p1);
\draw (p6) -- (p2);
\draw (p6) -- (p3);
\draw [thick](p6) -- (p4);
\end{tikzpicture}
\end{center}
\vspace{-0.3cm}
\caption{\emph{Flexible Bricard octahedron with point group $\mathcal{C}_s$.}}
\label{figurebricardcs}
\end{figure}

From the symmetry-extended version of Maxwell's rule, applied to the framework $(G,p)$ in Figure \ref{figurebricardcs}, the symmetry group $\mathcal{C}_s=\{Id,s\}$, and the homomorphism $\Phi_b:\mathcal{C}_s\to \textrm{Aut}(G)$ defined by
\begin{eqnarray}\Phi_b(Id)& = &id\nonumber\\ \Phi_b(s) &=& (v_{1}\,v_{3})(v_{2})(v_{4})(v_{5}\,v_{6})\textrm{,}\nonumber\end{eqnarray}
it follows that $(G,p)$, as well as any other $(\mathcal{C}_{s},\Phi_b)$-generic realization of $G$, has a fully $(\mathcal{C}_{s},\Phi_b)$-symmetric infinitesimal flex \cite{FGsymmax, BS2, BS4}. Thus, by  Corollary \ref{genflexcor}, any such realization of $G$ also has a $(\mathcal{C}_{s},\Phi_b)$-symmetry-preserving flex.

\begin{figure}[htp]
\begin{center}
\begin{tikzpicture}[very thick,scale=1]
\tikzstyle{every node}=[circle, fill=white, inner sep=0pt, minimum width=5pt];
\filldraw[fill=black!20!white, draw=black, thin, dashed]
    (-2.35,-0.45) -- (2.65,-1.49) -- (2.65,3.25) -- (-2.35,4.03) -- (-2.35,-0.45);
\filldraw[fill=black!20!white, draw=black, thin, dashed]
    (-0.7,-1.4) -- (0.7,-0.5) -- (0.7,4.1) -- (-0.7,3.2) -- (-0.7,-1.4);
\draw[dashed,thin](-2.35,-0.45) -- (2.65,-1.49);
\draw[dashed,thin](-2.35,4.03) -- (2.65,3.25);
\draw[dashed,thin](0.02,-0.95) -- (0.02,3.65);
\node [circle, shade, ball color=black!80!white, inner sep=0pt, minimum width=7pt](p1) at (-1.56,0.84) {};
\node [circle, shade, ball color=black!80!white, inner sep=0pt, minimum width=7pt](p3) at (1.96,0.086) {};
\node [circle, shade, ball color=black!80!white, inner sep=0pt, minimum width=7pt](p2) at (-0.4,-0.1) {};
\node [circle, shade, ball color=black!80!white, inner sep=0pt, minimum width=7pt](p4) at (0.4,0.41) {};
\node [circle, shade, ball color=black!80!white, inner sep=0pt, minimum width=7pt](p5) at (-1.14,2.85) {};
\node [circle, shade, ball color=black!80!white, inner sep=0pt, minimum width=7pt](p6) at (1.3,2.56) {};
\draw (p1) -- (-0.7,0.158);
\draw [dashed] (-0.7,0.158)--(p2)node[rectangle, fill=black!20!white, anchor=south, below=5pt] {$p_2$};
\draw [thick, dashed](p3) -- (p4) node[rectangle, fill=black!20!white,anchor=south, below=5pt] {$p_4$};
\draw (p2) -- (p3) node[rectangle, fill=black!20!white,anchor=south, right=5pt] {$p_3$};
\draw [thick, dashed](p4) -- (p1) node[rectangle, fill=black!20!white,anchor=south, left=5pt] {$p_1$};
\draw (p5)node[rectangle, fill=black!20!white,anchor=south, above=5pt] {$p_5$} -- (p1);
\draw (p5) -- (-0.7,1.18);
\draw [dashed](-0.7,1.18) -- (p2);
\draw (p5)-- (-0.7,2.48);
\draw [dashed](-0.7,2.48) -- (-0.01, 1.88);
\draw (-0.01, 1.88)-- (p3);
\draw [thick, dashed](p5) -- (p4);
\draw (p6)node[rectangle, fill=black!20!white,anchor=south, above=5pt] {$p_6$} -- (0,1.804);
\draw [dashed](0,1.804)--(-0.7,1.37);
\draw (-0.7,1.37)--(p1);
\draw (p6) -- (p2);
\draw (p6) -- (p3);
\draw [thick, dashed](p6) -- (p4);
\end{tikzpicture}
\end{center}
\vspace{-0.3cm}
\caption{\emph{Flexible octahedron with point group $\mathcal{C}_{2v}$.}}
\label{figurebricardc2v}
\end{figure}
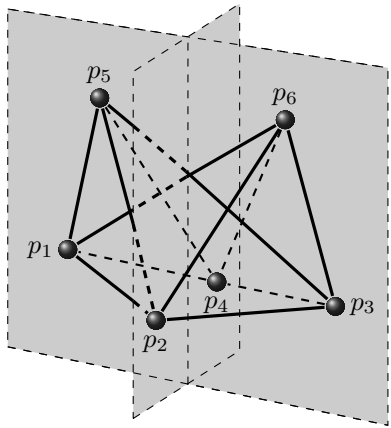

Finally, consider the $(\mathcal{C}_{2v},\Phi_c)$-generic framework $(G,p)$ in Figure \ref{figurebricardc2v}, where $\Phi_c:\mathcal{C}_{2v}\to \textrm{Aut}(G)$ is the unique type determined by the injective realization of $G$. Although $(G,p)$ is neither $(\mathcal{C}_{2},\Phi_a)$-generic nor $(\mathcal{C}_{s},\Phi_b)$-generic, we anticipate from the discussion above that $(G,p)$ possesses a flex that preserves both the $\mathcal{C}_{2}$ and the $\mathcal{C}_{s}$ symmetry defined in these examples.\\\indent The symmetry-extended version of Maxwell's rule applied to $(G,p)$,  $\mathcal{C}_{2v}$, and $\Phi_c$ detects a fully $(\mathcal{C}_{2v},\Phi_c)$-symmetric flex \cite{FGsymmax, BS2, BS4}.
Thus, by Corollary \ref{genflexcor}, the framework $(G,p)$, as well as any other $(\mathcal{C}_{2v},\Phi_c)$-generic realization of $G$, indeed possesses a $(\mathcal{C}_{2v},\Phi_c)$-symmetry-preserving flex.

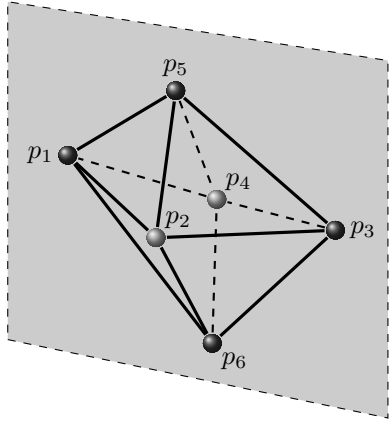
\begin{figure}[htp]
\begin{center}
\begin{tikzpicture}[very thick,scale=1]
\tikzstyle{every node}=[circle, fill=white, inner sep=0pt, minimum width=5pt];
\filldraw[fill=black!20!white, draw=black, thin, dashed]
    (-2.35,-0.45) -- (2.65,-1.49) -- (2.65,3.25) -- (-2.35,4.03) -- (-2.35,-0.45);
\node [circle, shade, ball color=black!80!white, inner sep=0pt, minimum width=7pt](p1) at (-1.56,1.994) {};
\node [circle, shade, ball color=black!80!white, inner sep=0pt, minimum width=7pt](p3) at (1.96,1) {};
\node [circle, shade, ball color=black!40!white, inner sep=0pt, minimum width=7pt](p2) at (-0.4,0.9) {};
\node [circle, shade, ball color=black!40!white, inner sep=0pt, minimum width=7pt](p4) at (0.4,1.41) {};
\node [circle, shade, ball color=black!80!white, inner sep=0pt, minimum width=7pt](p5) at (-0.14,2.85) {};
\node [circle, shade, ball color=black!80!white, inner sep=0pt, minimum width=7pt](p6) at (0.34,-0.5) {};
\draw (p1) --(p2)node[rectangle, fill=black!20!white, anchor=north, above right=3pt] {$p_2$};
\draw [thick, dashed](p3) -- (p4) node[rectangle, fill=black!20!white,anchor=south, above right=3pt] {$p_4$};
\draw (p2) -- (p3) node[rectangle, fill=black!20!white,anchor=south, right=5pt] {$p_3$};
\draw [thick, dashed](p4) -- (p1) node[rectangle, fill=black!20!white,anchor=south, left=5pt] {$p_1$};
\draw (p5)node[rectangle, fill=black!20!white,anchor=south, above=5pt] {$p_5$} -- (p1);
\draw (p5)  -- (p2);
\draw (p5)--  (p3);
\draw [thick, dashed](p5) -- (p4);
\draw (p6)node[rectangle, fill=black!20!white,anchor=south, below right=3pt] {$p_6$} -- (p1);
\draw (p6) -- (p2);
\draw (p6) -- (p3);
\draw [thick, dashed](p6) -- (p4);
\end{tikzpicture}
\end{center}
\vspace{-0.3cm}
\caption{\emph{An isostatic octahedron in $\mathscr{R}_{(G,\mathcal{C}_s,\Phi_d)}$.}}
\label{isobricocta}
\end{figure}

\begin{remark}
\emph{If $G$ is the graph of the octahedron, $\mathcal{C}_s$ is a symmetry group in dimension 3, and $\Phi_d:\mathcal{C}_s\to \textrm{Aut}(G)$ is defined by
\begin{eqnarray}\Phi_d(Id)& = &id\nonumber\\ \Phi_d(s) &=& (v_{2}\,v_{4})(v_{1})(v_{3})(v_{5})(v_{6})\textrm{,}\nonumber\end{eqnarray}
then $G$ is $(\mathcal{C}_{s},\Phi_d)$-generically isostatic. The framework $(G,p)$ in Figure \ref{isobricocta}, for example, is a realization of $G$ in $\mathscr{R}_{(G,\mathcal{C}_s,\Phi_d)}$ which is isostatic by Cauchy's Theorem \cite{Cauchy, dehn}.
}
\end{remark}

\begin{remark}\label{suspension}
\emph{The above rigidity analyses of symmetric octahedra can directly be extended to analyses of symmetric frameworks that consist of an arbitrary $2n$-gon and two `cone-vertices' that are linked to each of the joints of the $2n$-gon. These kinds of frameworks are also known as `double-suspensions', and are studied in \cite{consusp}, for example.
}
\end{remark}

A number of other interesting and famous examples of symmetric frameworks can also be proven to be flexible with the methods presented in this paper. These include Bottema's famous mechanism in the plane (see \cite{bottema}, for example), ring structures and reticulated cylinder structures in 3-space like the ones examined in \cite{FG4} and \cite{tarnaiely}, for example, and various types of bipartite frameworks in 3-space (see \cite{BS4}). Each of these structures possesses a `symmetry-preserving' flex, and it is precisely this kind of flex that our symmetry-based methods detect in each case. While detection of these flexes is not new, our new approach allows a much simpler verification of the flexes than  previous methods. New flexes can also be detected, and some will be presented in \cite{SchW}.

\section*{Acknowledgements}

We thank Walter Whiteley for numerous interesting and helpful discussions.

\end{document}